\documentclass[reqno,centertags]{amsart}
\usepackage{amsmath,amsthm,amscd,amssymb}
\usepackage{latexsym}
\usepackage{hyperref}
\usepackage{latexsym}
\usepackage{bbm}
\usepackage{color}

\newcommand{\UNIF}{\operatorname{UNIF}}
\DeclareMathOperator{\SC}{SC}
\DeclareMathOperator{\MP}{MP}
\DeclareMathOperator{\GW}{GW}
\DeclareMathOperator{\HP}{HP}
\newcommand{\KMK}{\operatorname{KMK}}
\DeclareMathOperator{\Pois}{Pois}

\DeclareMathOperator{\DVZ}{DVZ}
\newcommand{\Arcsine}{\operatorname{Arcsine}}

\def \d{{\tt d}}

\def \g{{\tt g}}
\def \c{{\tt c}}

\newcommand{\V}{\mathcal V}
\newcommand{\T}{\mathbb T}
\newcommand{\R}{\mathbb R}

\newcommand{\Sr}{\mathcal{S}}
\newcommand{\ii}{{\mathrm{i}}}
\renewcommand{\Re}{\operatorname{Re}}

\newcommand{\Sz}{\operatorname{Sz}}
\newcommand{\DG}{\operatorname{DG}}

\def \sur#1#2{\mathrel{\mathop{\kern 0pt#1}\limits^{#2}}}

\newcounter{smalllist}

\numberwithin{equation}{section}
\newtheorem{theorem}{Theorem}[section]
\newtheorem{proposition}[theorem]{Proposition}

\theoremstyle{definition}

\theoremstyle{remark}
\newtheorem{remark}{Remark}

\definecolor{Red}{rgb}{1,0,0}
\definecolor{Blue}{rgb}{0,0,1}

\begin{document}
\title
{On some gateways between sum rules}
\author{Fabrice Gamboa}
\address{ Universit\'e Paul Sabatier, Institut de Math\'ematiques de Toulouse,  31062-Toulouse Cedex 9, France and ANITI,}
\email{fabrice.gamboa@@math.univ-toulouse.fr}
\author{Jan Nagel}
\address{Technische Universit\"at Dortmund, Fakult\"at f\"ur Mathematik, 44227 Dortmund, Germany}
\email{jan.nagel@tu-dortmund.de}
\author{Alain Rouault}
\address{Laboratoire de Math\'ematiques de Versailles, UVSQ, CNRS, Universit\'e Paris-Saclay, 78035-Versailles Cedex France}
\email[corresponding author]{alain.rouault@uvsq.fr}

\keywords{Sum rules, Szeg{\H o} mapping, Verblunsky coefficients,  Delsarte-Genin mapping, Jacobi coefficients, relative entropy.
}

\subjclass{42C05, 47B36, 15B52, 34L05}

\date{\today}

\begin{abstract} 
We present correspondences induced by some classical mappings between measures on an interval and measures on the unit circle. More precisely, we link their sequences of orthogonal polynomial and their recursion coefficients. We also deduce some correspondences between particular equilibrium measures of random matrix ensembles. Additionally, we show that these mappings open up gateways between the sum rules associated with some classical models, leading to new formulations of several sum rules.
\end{abstract}


\maketitle
\section{Introduction}
The relation between orthogonal polynomials on the unit circle (OPUC) and orthogonal polynomials on the line (OPRL) is a longstanding problem. When a measure on the unit circle is mapped to a measure on the real line, what is the relation between the orthogonal polynomials related to these measures or their recursion coefficients? First results in this direction go back to Szeg{\H o}, who found a relation between the orthogonal polynomials when the mapping on the real line is the pushforward under $z\mapsto z+z^{-1}$, now called the Szeg{\H o} mapping, see \cite{simon2}, p. 880 for a historical account. The relation between the recursion coefficients was found by Geronimus: surprisingly, the so-called Verblunsky coefficients $(\alpha_k)_{k\geq 0}$ of the recursion on the unit circle appear in a decomposition of the Jacobi coefficients on the real line, forming an identity now known as the Geronimus relations. Since then, a variety of mappings have been studied, motivated from applications for orthogonal polynomials \cite{bessis1983orthogonality,costa2013orthogonal} operator theory \cite{golinskii2010guseinov,derevyagin2012cmv,cantero2016darboux} or signal processing \cite{delsarte1986split,delsartetri2,bracciali2005real}. 

Let us highlight the implications of such mappings and relations in particular on important identities in spectral theory called sum rules. Sum rules are identities between two nonnegative functionals of a probability measure $\mu$ compactly supported on $\mathbb R$ (resp. $\nu$ supported on $\mathbb T$). On the one hand, the first functional is an entropy-like functional with respect to some reference measure. On the other hand, the second functional is built from Jacobi coefficients (resp. Verblunsky coefficients) of $\mu$ (resp. $\nu$) and vanishes only for the reference measure. We call the left hand side (LHS) the spectral side and the right hand side (RHS) the coefficient side. 

The first historical example of such a sum rule is the classical Szeg{\H o}-Verblunsky identity, 
\begin{equation}
\label{SVsumrule}
\frac{1}{2\pi}\int_{0}^{2\pi}\log g_{\nu}(\theta)d\theta = \sum_{k= 0}^\infty \log(1-|\alpha_k|^2)\,,
\end{equation}
where $\nu$ is a measure on the unit circle with Lebesgue decomposition 
\[d\nu(\theta)= g_{\nu}(\theta) \tfrac{d\theta}{2\pi}+d\nu_s(\theta)\]
having Verblunsky coefficients $(\alpha_k)_{k\geq 0}$. Both sides of \eqref{SVsumrule} vanish if, and only if, $\nu$ is the uniform measure on the circle (the reference measure in this case). We refer to Chapter 1 of \cite{simon2} for a discussion of the origin of this sum rule. 

The most famous sum rule for measures on the line is the Killip-Simon sum rule \cite{killip2003sum}. 
An exhaustive discussion and history of this sum rule can be found in  Section 1.10 of the book  \cite{simon2} and a deep analytical proof in Chapter 3 therein. The reference measure for this sum rule is the semicircle law (SC). 

An important consequence of these two sum rules is the equivalence of two conditions for the finiteness of both sides, one formulated in terms of Verblunsky or Jacobi coefficients and the other as a spectral condition. In the words of Simon \cite{simon2}, these are the \emph{gems} of spectral theory.
In \cite{gamboacanonical} and \cite{GaNaRo}, we revisited these results from a probabilistic point of view and gave a new proof based on large deviations. We also refer to the work of Breuer et al. \cite{BSZ} which enlightens non-probabilists about  \cite{GaNaRo}, \cite{gamboacanonical}. The method was robust enough to prove new sum rules with reference measures such as Marchenko-Pastur (MP), Kesten-McKay (KMK) on the real line and Gross-Witten (GW), Hua-Pickrell (HP) on the unit circle.

The main contribution of this paper is two-fold. On the one hand, we gather a series of results on relations between measures on the unit circle and measures on the real line and their orthogonal polynomials under several well known  mappings: Szeg{\H o}, Delsarte-Genin (DG), Derevyagin-Vinet-Zhedanov (DVZ) and M\"obius. On the other hand, we show how these relations allow to catch a --potentially new-- sum rule from an existing one. The main idea is easy: we transform  both sides  according to the mapping. While our first contribution is merely expository in nature, we believe the second contribution can be of great interest, either to find new sum rules or to highlight connections, or ``gateways'', between existing identities. 

As an easy example for such a gateway, the Szeg{\H o}-Verblunsky sum rule \eqref{SVsumrule} leads to an identity for measures on $[-2,2]$, when both sides are transformed according to the Szeg{\H o} mapping. The LHS may be written as an integral with respect to the Arcsine law while the Geronimus relations allow to rewrite the RHS (see Section \ref{susec:gatewayUNIFarcsine}).   
To give an overview of further results (we refer to Section 3 for the statement of the sum rules and Section 4 for the definition of the mappings): 
\begin{itemize}
\item
Particular cases of the KMK-sum rule can be obtained from the HP-sum rule by the Szeg{\H o} mapping or by the DG mapping (Section \ref{susec:gatewayHPKMK}).
\item
The GW-sum rule implies the new sum rules \eqref{newSRGW} and \eqref{newSRGW-1} by the Szeg{\H o} mapping. 
\item
The GW-sum rule leads to the reformulations \eqref{nsr2} and \eqref{newsr2b} under the DG mapping, with a new formula for the RHS in Theorem \ref{newrhog}.
\item
Under the DVZ mapping, the GW-sum rule leads to a variant of the Killip-Simon sum rule \eqref{KSSRvariant}. 
\item
We prove a new sum rule with reference to the Poisson measure $\Pois(\zeta)$ in Theorem \ref{Pois}.
\item
Another new Poisson sum rule is obtained from a recent result of \cite{bessonov} in Proposition \ref{rembess}. 
\item
A new analytical proof of a weak version of the HP-sum rule is in Proposition \ref{weakprop}.
\end{itemize}

The  paper is organized as follows: In Section 2 we recap the required background on orthogonal polynomials on the real line and on the unit circle with corresponding recursion formulas. In Section 3 we recall the main known  sum rules with their reference measures. Section 4 discusses the  mappings used in our work. Then,  in Section 5 we apply these mappings to our reference measures. Section 6 presents the gateways between OPUC sum rules and OPRL sum rules. In Section 7 are the proofs of some new sum rules and some auxiliary results are in Section 8.


\section{Orthogonal polynomials}
Let $\mathcal{M}_1(\mathbb{R})$ (resp. $\mathcal{M}_1(\mathbb{T})$) denote the set of all probability measures on $\mathbb{R}$ (resp. on the unit circle $\mathbb{T} = \partial \mathbb{D}$, where $\mathbb{D}$ is the open unit disk $\mathbb{D} = \{z\in \mathbb{C}: |z|< 1\}$). Additionally, we write $\mathcal{M}_{1,s}(\mathbb{T})$ for the set of all symmetric probability measures on $\mathbb{T}$, invariant under the transformation $z\mapsto \bar z$.

\subsection{OPRL}

The sequence of orthogonal polynomials on the real line (OPRL) is well defined for a probability measure $\mu\in \mathcal{M}_1(\mathbb{R})$ with a compact support consisting of infinitely many points, a.k.a. nontrivial case  (in constrast to a finite support consisting of $n$ points, a.k.a. trivial case). They are obtained by applying
to the sequence $1, x, x^2, \dots$ the orthonormalizing Gram-Schmidt procedure. The resulting polynomials $p_0,p_1,\dots$, with $p_k$ of degree $k$, obey the recursion relation
\begin{align} \label{polrecursion}
xp_k(x) = a_{k+1} p_{k+1}(x) + b_{k+1} p_k (x) + a_{k} p_{k-1}(x)
\end{align}
for $ k \geq 0$, with $p_{-1}=0$. The recursion or Jacobi coefficients (or short ``J-coefficients'') satisfy that for all $k$, $b_k \in \mathbb R$ and  $a_k > 0$. Notice that here the orthogonal polynomials are not monic but normalized in $L^2(\mu)$. The monic polynomials satisfy the recursion
\begin{align} \label{polrecursionmonic}
xP_k(x) = P_{k+1}(x) + b_{k+1} P_k (x) + a_{k}^2 P_{k-1}(x)\,.
\end{align}
When the support of $\mu$ consists of $n$ points, the orthogonal  polynomials  might be defined up to degree $n-1$ and J-coefficients $b_1,a_1,\dots ,a_{n-1},b_n$ are well defined.

For a non-trivial measure $\mu$ let us equip the vector space $L^2(\mu)$ with the 
basis $(p_k)_{k\geq 0}$.
 Then the linear map $f \mapsto xf$, multiplication by the identity, is represented by the tridiagonal matrix  
\begin{align}
\label{favardinfini}
J_\mu = \begin{pmatrix} b_1&a_1 &0&0&\cdots\\
a_1&b_2  &a_2&0&\cdots\\
0&a_2 &b_3&a_3& \\
\vdots& &\ddots&\ddots&\ddots
\end{pmatrix} .
\end{align}
Conversely, if $H$ is a bounded Hermitian operator on an infinitely dimensional Hilbert space $\mathcal H$, and $e$ is a cyclic vector, then we can define the spectral measure $\mu$ of the pair $(H, e)$ and
then $(\mathcal H, H,e)$ is isomorphic to $(\ell^2, J_\mu, e_1)$ where $e_1= (1, 0, 0 , \dots)^t$. 
Such a correspondence still holds between Hermitian operators on an $n$-dimensional space and measures supported by $n$ points and $n\times n$ tridiagonal matrices.   

If the support of $\mu$ is contained in $[0, \infty)$, there is a decomposition of J-coefficients,  
\begin{align}
\label{zcoeff}
\notag
b_k &= z_{2k-2} + z_{2k-1},\\
a_k^2 &= z_{2k-1}z_{2k} ,
\end{align}
with $z_0 =0$ and $z_k \geq 0$ for $k\geq 1$ (see \cite{Chihara} p.47).
The $z_k$ will be called the canonical coefficients and they are uniquely determined by the J-coefficients.

If $\mu$ is nontrivial with support contained in the interval $[-2,2]$, there exists a decomposition of J-coefficients,
\begin{align}\label{eq:coeffdecomposition}
\notag
b_{k+1} &= (1-u_{2k})u_{2k+1} -  (1+u_{2k})u_{2k-1},\\
a_{k+1}^2 &= (1-u_{2k})(1-u_{2k+1}^2)(1+u_{2k+2}) ,
\end{align}
with $u_{0}=-1$ and $u_k\in (-1,1)$ for all $k\geq 1$. The $u_k$ will be called canonical moments (although more classically, the $\tfrac{1}{2}(u_k+1)$ are called canonical moments \cite{DeSt97}) and they are uniquely determined by the  J-coefficients. If $\mu$ is nontrivial,  $u_k\in (-1,1)$ for all $k\geq 1$, while if $\mu$ is supported by $n$ points, we can still define $u_1,\dots ,u_{2n-2}\in (-1,1)$ and $u_{2n-1}\in \{-1,1\}$. Let us notice that if $\mu$ is symmetric, then $u_{2k+1} =0$ for all $k$, all  the diagonal coefficients $b_k$ vanish and
\begin{align}\label{ccUNIF}a_{k+1}^2 = (1 - u_{2k}) (1 + u_{2k+2})  \ \  (k \geq 0)\,.\end{align}

\subsection{OPUC}
For a probability measure $\nu\in \mathcal{M}_1(\mathbb{T})$  supported by at least $k+1$ points, the inductive relation  between two successive monic polynomials $\Phi_{k+1}$ and $\Phi_{k}$, where $\Phi_{k}$ has degree $k$, orthogonal with respect to $\nu$ 
involves a complex number $\alpha_k$
and may be written as
\begin{equation}
\label{recpolycirc}
\Phi_{k+1}(z)=z\Phi_{k}(z)-\overline{\alpha}_k\Phi_{k}^*(z)\mbox{ where } \Phi_{k}^*(z):=z^k\overline{\Phi_k(1/\bar{z})}.
\end{equation}
The complex numbers $\alpha_k=-\overline{\Phi_{k+1}(0)} \ , k \geq 0$ are the so-called Verblunsky coefficients (in short V-coefficients). They are  also called Schur, Levinson, Szeg\H{o} coefficients in other contexts or canonical moment as well \cite{DeSt97}. 
We also set $\alpha_{-1}=-1$.
The V-coefficients satisfy $|\alpha_{k-1}|<1$ if $k\geq 1$ and the support of $\nu$ contains at least $k+1$ points and $|\alpha_{k-1}|=1$ if the support consists of exactly $k$ points. For a symmetric measure $\nu \in \mathcal{M}_{1,s}(\mathbb{T})$, the V-coefficients are real. We will denote by $(\varphi_k)_{k\geq 0}$ the sequence of orthonormal polynomials on the unit circle (OPUC).

In the basis $(\chi_k)_{k \geq 0}$ obtained by orthonormalizing $1, z , z^{-1}, z^2, z^{-2}, \dots$,
the linear transformation $f \rightarrow zf$ in $L^2(\nu)$ is represented by the so-called CMV-matrix
\begin{align}
\label{favardfini}
\mathcal C_\mu = \begin{pmatrix} \bar \alpha_0&\bar \alpha_1\rho_0 &\rho_1\rho_0&0&0&\dots\\
\rho_0& -\bar\alpha_1\alpha_0&-\rho_1\alpha_0 &0 &0 &\dots\\
0&\bar\alpha_2\rho_1&-\bar\alpha_2\alpha_1&\bar\alpha_3\rho_2&\rho_3\rho_2&\dots\\
0&\rho_2\rho_1& -\rho_2\alpha_1& -\bar\alpha_3\alpha_2&-\rho_3\alpha_2&\dots\\
0&0&0&\bar\alpha_4\rho_3&-\bar\alpha_4\alpha_3& \dots\\
\dots&\dots&\dots&\dots&\dots&\dots
\end{pmatrix} 
\end{align}
with $\rho_k = \sqrt{1- |\alpha_k|^2}$ for every $k \geq 0$ in the non-trivial case.

Conversely, if $U$ is a unitary operator on an infinite dimensional Hilbert space $\mathcal H$ and $e$ is a cyclic vector, then we can define the spectral measure $\nu$ of the pair $(U, e)$ and
then $(\mathcal H, U,e)$ is isomorphic to $(\ell^2, \mathcal C_\nu, e_1)$. 
%
%
%
%
Let
\begin{align}
\Theta_k = \begin{pmatrix} \alpha_k & \rho_k\\
\rho_k& - \alpha_k
\end{pmatrix}
\end{align}
and
\begin{align}
\label{defLM}
\mathcal L = \Theta_0\oplus\Theta_2\oplus \cdots ,\qquad  \mathcal M = {\bf 1}\oplus \Theta_1\oplus\Theta_3 \oplus \cdots ,
\end{align}
where ${\bf 1}$ denotes the $1\times 1$ identity matrix and $\oplus$ is the direct sum operator. The unitary operators $\mathcal L$ and $\mathcal M$ satisfy
\begin{align} \label{CMVdecomposition}
\mathcal C_\nu = \mathcal L \mathcal M\,.
\end{align}

For probability measures $\nu,\mu$ both on $\mathbb{R}$ or on $\mathbb{T}$, let $\mathcal K (\nu | \mu)$ denote the Kullback-Leibler divergence or relative entropy of $\nu$ with respect to $\mu$: 
\begin{equation}
\label{KL}
{\mathcal K}(\nu\, |\, \mu)= \begin{cases}  \ \displaystyle\int
\log\frac{d\nu}{d\mu}\!\ d\nu\;\;& \mbox{if}\ \nu\ll\mu  \ \hbox{ and }\,
\log \frac{d\nu}{d\mu}\in L^1(\nu), \\
   \  \infty  &  \mbox{otherwise.}
\end{cases}
\end{equation}

\section{Reference measures and sum rules}
\label{sec:measuressumrules}

\subsection{Measures on $\mathbb{R}$}

We start with measures on the real line and state sum rules relative to these measures. In order to formulate the spectral side, we need some support conditions. 
For $c^-<c^+$ we define the set $\mathcal{S}_1(c^-,c^+)$ as the set of probability measures $\mu$ on $\mathbb{R}$ whose support satisfies 
\begin{align*}
\operatorname{supp}(\mu) = I \cup E ,
\end{align*}
where $I\subset [c^-,c^+]$ and $E=E(\mu)$ is an at most countable subset of $[c^-,c^+]^c$.

\subsubsection{Semicircle distribution}

The semicircle law is 
\begin{align}
\SC(dx) = \frac{1}{2\pi} \sqrt{4 -x^2}\, \mathbbm{1}_{\{-2 \leq x \leq 2\} }\ dx .
\end{align}
It is the central probability measure in classical random matrix theory. 
Indeed, it is the equilibrium measure for
a large class of random matrix models  
(the limit of their empirical eigenvalue distribution). The Jacobi matrix of $\SC$ is also called the free Jacobi matrix with J-coefficients 
\begin{equation}
\label{JcH}a^{\SC}_k = 1,\  b^{\SC}_k =0 \ \ \hbox{for all} \ k\geq 1\,.
\end{equation}
We start by stating the classical sum rule of \cite{Killip2} (and explained in \cite{Simon-newbook} p.37), the new probabilistic proof using large deviations might be found   in \cite{GaNaRo}.  
For a probability measure $\mu$ on $\R$ with recursion coefficients $a_k,b_k$ as in \eqref{polrecursion}, define the sum
\begin{align}\label{rateG1}
\mathcal{I}_H(\mu) = \frac{1}{2}\sum_{k= 1}^\infty  b_k^2 + G(a_k^2) ,
\end{align}
where  
\begin{align}
\label{defG}
G(x) = x -1 - \log x\,.
\end{align} 
Furthermore, define 
\begin{align*}
\mathcal{F}_{\SC}(x) :=  \displaystyle \int_2^{|x|} \sqrt{t^2-4}\!\ dt = \displaystyle\frac{|x|}{2} \sqrt{x^2-4} - 2 \log \left(\displaystyle \frac{|x|+\sqrt{x^2-4}}{2}\right)
\end{align*}
if $|x|\geq 2 $ and $\mathcal{F}_H(x) = \infty$ otherwise. Then the following remarkable identity holds.

\medskip

\begin{theorem}[\cite{Killip2}]
\label{sumruleg}
Let $J$ be a Jacobi matrix with diagonal entries $b_1,b_2,\ldots \in \R$ and subdiagonal entries $a_1,a_2,\ldots >0$ satisfying  $\sup_k (a_k + |b_k|) < \infty$ and let $\mu$ be the associated spectral measure. Then $\mathcal{I}_H(\mu)$ is infinite if $\mu \notin \Sr_1(-2,2)$ and for $\mu \in \Sr_1(-2,2)$, 
\begin{align*} 
 {\mathcal K}(\operatorname{SC}\, |\!\ \mu) +  \sum_{\lambda \in E(\mu)} {\mathcal F}_{\SC}(\lambda) = \mathcal{I}_H(\mu), 
\end{align*}
where both sides may be infinite simultaneously. 
\end{theorem}

Let us emphazise that for a sum rule as in Theorem \ref{sumruleg}, both sides are nonnegative and vanish if and only if $\mu$ is equal to the reference measure, which is the semicircle law $\operatorname{SC}$ in this case.

\subsubsection{Marchenko-Pastur distribution}

The Marchenko-Pastur distribution with parameter $\tau \in (0,1]$ is
\begin{align}
\MP_\tau (dx) = \frac{\sqrt{(x-\tau_-)(\tau_+ -x)}}{2\pi\tau x} \mathbbm{1}_{[\tau_-, \tau_+]}(x) dx\,.
\end{align}
where $\tau_\pm = \sqrt{1\pm\tau}$. In random matrix theory, it is the equilibrium measure of the Laguerre ensemble. 
Its canonical coefficients (see \eqref{zcoeff}) are
\[z_{2k-1} = 1  , \quad z_{2k} = \tau \ (k \geq 1) \]
with $z_0=0$, which correspond to the J-coefficients
\[a_k = \sqrt \tau  ,\quad  b_{k+1} = 1+\tau \ (k\geq 1)\]
with $b_1=1$. 
Notice that the MP distribution is not symmetric.
Let $\mathcal I_L$ be defined by
\begin{align}
\mathcal I_L (\mu) = \sum_{k=1}^\infty \tau^{-1}G(z_{2k-1}) + G(\tau^{-1}z_{2k})\,.
\end{align}
Furthermore, define for $x \notin (\tau^-, \tau^+)$ 
\begin{align*}
\mathcal{F}_{\MP}(x) := \int_{I(x)} \frac{\sqrt{(t-\tau^+)(t-\tau^-)}}{\tau t}\!\ dt
\end{align*}
where $I( x) = [\tau^+, x]$ if $x\geq \tau^+$ and $I(x) = [x, \tau^-]$ if $0 < x \leq \tau^-$. 
 Then we have the following theorem.
\begin{theorem}[\cite{GaNaRo} Theorem 2.2]
\label{sumrulelag}
Let $\mu \in \mathcal M_1([0, \infty))$ be a non-trivial measure with compact support and $0 < \tau \leq 1$.  Then $\mathcal I_L(\mu) = \infty$ if  $\mu \notin  \Sr_1(\tau^-, \tau^+)$ and if $\mu \in \Sr_1(\tau^-, \tau^+)$ we have
\begin{align}
\label{SRL} 
{\mathcal K}(\operatorname{MP(\tau)} | \mu) +  \sum_{\lambda \in E(\mu)} {\mathcal F}_{\MP}(\lambda) = \mathcal{I}_L(\mu) , 
\end{align}
where both sides may be infinite simultaneously. 
\end{theorem}

\subsubsection{The KMK distribution}

The Kesten-McKay law with parameters $\kappa_1,\kappa_2\geq 0$ is denoted by $\KMK(\kappa_1, \kappa_2)$ and has the density
\begin{align}
\KMK(\kappa_1, \kappa_2)(dx) = \frac{(2+\kappa_1+\kappa_2)}{2\pi} \frac{\sqrt{(x-u_-)(u_+-x)}}{4-x^2}\mathbbm{1}_{\{u^- <x < u^+\}}\ dx 
\end{align}
where
\begin{align}
u^\pm = \frac{2\left(\kappa_2^2 - \kappa_1^2\pm 4\sqrt{(1+\kappa_1)(1+\kappa_2)(1+\kappa_1+\kappa_2)}\right)}{(2+ \kappa_1+\kappa_2)^2}\,.
\end{align}
It is the equilibrium measure of the Jacobi ensemble. 
The canonical moments of $\KMK(\kappa_1,\kappa_2)$ of even and odd index are, respectively:
\begin{align}
\label{KMKcanonical12}
u_{2k}^{\kappa_1, \kappa_2} \equiv u_e^{\kappa_1, \kappa_2} := -\frac{\kappa_1 + \kappa_2}{2+ \kappa_1+\kappa_2} \ , \ u_{2k+1}^{\kappa_1, \kappa_2} \equiv u_o^{\kappa_1, \kappa_2} := \frac{\kappa_2-\kappa_1}{2+ \kappa_1+\kappa_2}\,.
\end{align}
We will consider also a symmetric version with $\kappa_1= \kappa_2=\kappa$, we will denote it $\KMK(\kappa) := \KMK(\kappa, \kappa)$:
\begin{align}  \label{defKMK}
\KMK(\kappa)(dx) = \frac{(1+\kappa)}{\pi}\frac{\sqrt{u^2 -x^2}}{4-x^2} \mathbbm{1}_{\{|x| \leq u\}}\ dx ,
\end{align}
where
\[
 u=2\frac{\sqrt{1+2\kappa}}{1+\kappa}\,.\]
The canonical moments of $\KMK(\kappa)$ 
 are  
\begin{align}\label{KMKcanonical}
u^{\kappa}_{2k} =u^{(\kappa)} :=\frac{-\kappa}{1+\kappa} ,\quad   u^{(\kappa)}_{2k-1} = 0\,,
\end{align}
see \cite[Section 6]{gamboa2011large} for the linearly transformed canonical moments. 
For $\kappa=0 $, the Kesten-McKay law is the $\Arcsine$ distribution
\begin{align}
\label{defAS}
\Arcsine(dx)  := \frac{1}{\pi \sqrt{4-x^2}} \mathbbm{1}_{\{-2 < x < 2\}}\ dx \,,
\end{align}
whose canonical coefficients are all zero.

To state the sum rule, we need  some more
notation. Set for  $u\in (-1,1)$
\begin{align} \label{Hoddeven}
\begin{split}
\mathcal H_e^{\kappa_1, \kappa_2}(u) &:=  - (1+\kappa_1 + \kappa_2)\log\frac{1-u}{1-u_e^{\kappa_1, \kappa_2}} - \log \frac{1+u}{1+u_e^{\kappa_1, \kappa_2}}
 \,, \\
 \mathcal H_o^{\kappa_1, \kappa_2} (u) &:= -(1+ \kappa_1)\log\frac{1-u}{1-u_o^{\kappa_1, \kappa_2}} -(1+ \kappa_2)\log \frac{1+u}{1+u_o^{\kappa_1, \kappa_2}}
 \,.
\end{split}
\end{align}

For a non-trivial measure $\mu \in \mathcal{M}_1([-2,2])$ with canonical moments $u_k\in (-1,1)$, define 
\begin{align} \label{rateJacobicoefficient}
\mathcal{I}_J(\mu) =\sum_{k=1}^\infty \mathcal H_o^{\kappa_1, \kappa_2} ( u_{2k-1}) +  \mathcal H_e^{\kappa_1, \kappa_2} (u_{2k})\,.
\end{align} 
Finally, for the contribution of the outlying support points, we define for $x \notin (u^-, u^+)$
\begin{align} \label{outlierJac}
{\mathcal F}_{\KMK(\kappa_1, \kappa_2)}(x) = \int_{I(x)} (2+\kappa_1+\kappa_2)\frac{\sqrt{(t-u^+)(t-u^-)}}{4-t^2}\!\ dt
\end{align}
where $I(x) = [u^+, x]$ if $x \in [u^+, 2]$ and $I(x) = [x, u^-]$ if $x \in [-2, u^-]$.

We are now able to give the sum rule relative to the KMK measure. It is Theorem 2.3 in \cite{GaNaRo}, where it is formulated for linearly transformed measures on $[0,1]$.

\begin{theorem}
\label{LDPSR}
Let $\mu \in \mathcal{M}_1([-2,2])$ be a nontrivial measure and $\kappa_1, \kappa_2 \geq 0$. Then $\mathcal{I}_J(\mu) =\infty$ if $\mu \notin \mathcal S_{1}(u^-,u^+)$, and if $\mu\in \mathcal{S}_1(u^-,u^+)$ we have
\begin{align} 
\label{sumrule}
\mathcal K(\KMK(\kappa_1, \kappa_2)\, |\, \mu) + \sum_{\lambda \in E(\mu)} \mathcal F_{\KMK(\kappa_1, \kappa_2)} (\lambda) = \mathcal{I}_J(\mu) ,
\end{align} 
where both sides may be infinite simultaneously. 
\end{theorem}

\begin{remark}
In the particular case $\kappa_1=\kappa_2 =0$, we obtain 
 the identity
\begin{align}
\label{SRarcsine}
\mathcal K( \Arcsine\, | \, \mu) = - \sum_{k=1}^\infty \log (1-u_k^2)\,.
\end{align}
It is very close to the results of Gamboa and Lozada \cite{gamboa2004large} and equivalent to the so-called $C_0$ sum rule of Simon and Zlatos \cite[Theorem 13.8.8]{simon2}:
\begin{align}
\label{Zlatos}
\mathcal K( \Arcsine\, | \,\mu) = \log 2 - \sum_{k=1}^\infty \log a_k^2\,. 
\end{align}
Indeed, using \eqref{eq:coeffdecomposition} we can write
\begin{align}
\label{2series}
-\sum_{k=1}^n \log a_k^2 
= -\log 2  - \log(1 - u_{2n}) + \sum_{k=1}^{2n}  -\log (1 - u_k^2)\,.
\end{align} 
Suppose that the last sum is bounded, then $\lim_{n\to \infty}  u_{2n} =0$ and hence
\begin{align}
\label{2seriesa}
-\sum_{k=1}^\infty \log a_k^2 =\log 2 -\sum_{k=1}^\infty \log(1 -u_k^2)\,. 
\end{align}
On the other hand, since $\log (1 - u_{2n+2}) \leq \log 2$ , \eqref{2series} implies
\[\sum_{k=1}^{2n}  -\log (1 - u_k^2)\leq 2 \log 2 -\sum_{k=1}^n \log a_k^2\]
and when the sum on the LHS diverges, $-\sum_{k=1}^\infty \log a_k^2$ does as well, so that \eqref{2seriesa} holds true in any case.
\end{remark}

\subsection{Measures on $\mathbb{T}$}

In analogy to the real case we introduce for $0\leq \theta^-<\theta^+\leq 2\pi$ the set $\mathcal{S}^{\mathbb{T}}_1(\theta^-,\theta^+)$ of probability measures $\nu \in \mathcal{M}_1(\mathbb{T})$ supported on $I\cup E$, where $I$ is a subset of the arc
\begin{align} \label{circlearc}
\{z=e^{\ii \theta} \in \mathbb{T}\!\ |\, \theta \in [\theta^-,\theta^+]\}
\end{align}
and where $E=E(\nu)$ is an at most countable subset of the complement of the set \eqref{circlearc}. 

\subsubsection{Uniform distribution}

We write $\UNIF$ for the normalized Lebesgue measure on $\mathbb T$
\[\UNIF (d\theta) = \frac{d\theta}{2\pi}\,.\]
Its V-coeffcients are
\[\alpha_k = 0 , \quad k\geq 0\,.\]
The classical Szeg\H{o}-Verblunsky theorem (see \cite{Simon-newbook}, Theorem 1.8.6) is the identity
\begin{equation}
\label{segverth0}
\frac{1}{2\pi} \int_{\mathbb{T}} \log g_{\nu}(e^{\ii \theta})d\theta = \sum_{k=0}^\infty\log(1-|\alpha_k|^2)\,,
\end{equation}
where $\nu\in \mathcal{M}_1(\mathbb{T})$ is nontrivial with V-coefficients $\alpha_k$ and with Lebesgue decomposition 
\begin{align*}
d\nu= g_{\nu} d\UNIF +d\nu_s
\end{align*}
with respect to $\UNIF$. Changing  signs in both sides of this equation leads to 
\begin{equation}
\label{SVsum}
\mathcal K(\UNIF | \nu) = -\sum_{k=0}^\infty \log(1 - |\alpha_k|^2) .
\end{equation}

\subsubsection{Gross-Witten} 
\label{GWdefi}

The Gross-Witten measures are a class of equilibrium measures for random matrix distributions with potential  
\begin{equation}
\label{potGW}
\V_\g(z) =  - \g \Re (z) \,,
\end{equation}
with parameter $\g \in \mathbb{R}$. 
For details and applications of this distribution we refer to \cite{HiaiP} p. 203, \cite{GrossW}, and \cite{Wadia}. 

If $-1\leq \g \leq 1$ (ungapped or strongly coupled phase), the Gross-Witten measure $\GW(\g)$ is supported by $\mathbb T$ and is given by :
\begin{equation}
\label{GW-}
\GW(\g)(dz)
= \frac{1}{2\pi} (1 + \g \cos \theta)\!\ d\theta, 
\end{equation}
with $z = e^{\ii \theta }, \theta \in [-\pi, \pi )$. Note that $\tau_\pi(\GW(\g)) = \GW(-\g)$, 
where
\begin{align}
\label{rot}
\int f(\theta) d\tau_\pi (\mu)(\theta) = \int f(\theta +\pi) d\mu(\theta)\,.\end{align}
Since\footnote{In the sequel, we will use the notation  $\alpha_k(\nu)$ or $u_k(\mu)$ when the context needs the name of the measure we work with.} 
\begin{align} \label{GW+-}
\alpha_k(\tau_\pi (\mu)) = (-1)^{k+1}\alpha_k (\mu) ,
\end{align}
see \cite{simon05}, we state the V-coefficients only for the case $\g<0$.

For $-1 \leq \g < 0$, the measure $\GW(\g)$ has V-coefficients
\begin{equation}
\label{alphalimGW-}
\alpha_n^\g = \alpha_n(\GW(\g)) = \begin{cases}\displaystyle  -\frac{x_+ - x_-}{x_+^{n+2} - x_-^{n+2}} & \hbox{if} \ 
-1 <\g< 0, \\
\displaystyle -\frac{1}{n+2}& \hbox{if} \ \g = -1\,,
\end{cases}
\end{equation}
(see Simon \cite{simon05}, p. 86),
where $x_\pm = -\g^{-1} \pm \sqrt{\g^{-2}-1}$
 are roots of the equation 
\[x + \frac{1}{x} = -\frac{2}{\g}\,.\]
We remark that the measure $\GW(\g)$ has only nontrivial moments of order  $\pm 1$.

For $|\g| \geq 1$  (gapped or weakly coupled phase), let $\theta_{g}  \in [0, \pi]$ be the solution of
\begin{equation}\label{eqthetag}
\sin^2 (\theta_\g /2) = |\g|^{-1}\,.
\end{equation}
When  $\g \leq - 1$,  the
Gross-Witten measure is  for $z = e^{\ii \theta} , \theta \in [0, 2\pi)$
\begin{equation}\label{GWMeq}
\GW(\g) (dz) = \frac{|\g|}{\pi}\sin(\theta/2)\!\  \sqrt{\sin^2(\theta/2)-\cos^2(\theta_\g/2)}\!\ 1_{[\pi-\theta_\g, \pi+\theta_\g]}\!\ d\theta\,. 
\end{equation}
Summarizing formula (7.22) in Zhedanov \cite{Zhedanov}, we have that in the case $\g<-1$ the V-coefficients are 
\begin{align}\label{alphalimGW+}
\alpha_{n-1}^\g = \alpha_{n-1}(\GW(\g)) 
= 1 - \frac{2}{1+q}\frac{1 - q^{n+2}}{1 - q^{n+1}}\,,
\end{align}
where
\begin{align}\label{GWlim}
q &= \left(\sqrt{|\g|} - \sqrt{|\g| -1}\right)^2\, .
\end{align}
Since $0 < q<1$, it holds that
\begin{align}
\label{limangle}
\lim_{n\to \infty} \alpha_n^\g = - \sqrt{1- |\g|^{-1}} = -\cos(\theta_\g/2)\,.
\end{align}
When $g\geq 1$, the equilibrium measure is 
\begin{align}
\GW(\g) (dz) = \frac{|\g|}{\pi}\cos(\theta/2)\!\  \sqrt{\sin^2(\theta_\g/2)-\sin^2(\theta/2)}\!\ 1_{[-\theta_\g, \theta_\g]}\!\ d\theta\,. 
\end{align}
Note that again $\tau_\pi(\GW(\g)) = \GW(-\g)$, so that by \eqref{GW+-} the V-coefficients in this case can be obtained from \eqref{alphalimGW+}.

\begin{remark}
We may rotate GW as in \cite{ARH} and consider the equilibrium measure obtained by pushing forward $\GW(\g)$ by a rotation of angle $\eta$ instead of $\pi$ in \eqref{rot}. 
\end{remark}

The first sum rule relative to the Gross-Witten equilibrium measure $\GW(\g)$ was discovered by Simon for $\g = -1$ (see  \cite [Theorem 2.8.1]{simon05}), proved later with probabilistic methods by Breuer, Simon and Zeitouni \cite{breuer2018large}.  It is easily extended to $ |\g| \leq 1$ ( \cite[Corollary 5.4]{GNROPUC}).
For $-1 \leq \g \leq 0$ and  $\nu \in \mathcal{M}_1(\mathbb{T})$ nontrivial, it is the identity
\begin{align}
\notag
\mathcal K(\GW(\g)\ | \ \nu) &= H(\g)   +\frac{\g }{2}-\frac{\g }{2} \sum_{k=0}^\infty|\alpha_k - \alpha_{k-1}|^2  \\
\label{GWclassical}
&\quad + \sum_{k=0}^\infty - \log (1 - |\alpha_k|^2)  + \g|\alpha_k|^2 \,,
\end{align}
where
\begin{align}
\label{Hvalue}
H(\g) = \mathcal K(\GW(\g) \ | \ \UNIF) = 1- \sqrt{1- \g^2} + \log \frac{1 + \sqrt{1-\g^2}}{2} .
\end{align}
We recall that in \eqref{GWclassical}, $\alpha_{-1}=-1$. 
The sum rule \eqref{GWclassical} implies the following \emph{gem}, conditions for finiteness of the Kullback-Leibler divergence. 
The RHS of \eqref{GWclassical} is finite if and only if
\begin{align}
\label{gem2}
\sum_{k=0}^\infty |\alpha_k|^2 < \infty \ &\hbox{if} \ -1 <\g \leq 0, \\
\label{gem4}
\sum_{k=0}^\infty |\alpha_k|^4 < \infty \ &\hbox{and} \ \sum_{k=1}^\infty |\alpha_k- \alpha_{k-1}|^2 < \infty \ \hbox{if}\  \g =-1\, .
\end{align}

\begin{remark}
Provided that $\sum_k |\alpha_k|^2<\infty$, we may rewrite the sum rule as
\begin{align}
\mathcal K(\GW(\g) \ | \ \nu) =  H(\g)+ \g \Re \sum_{k=0}^\infty \alpha_k \bar\alpha_{k-1} - \sum_{k=0}^\infty \log(1- |\alpha_k|^2) ,
\end{align}
 \cite[p. 174]{simon05},  for the case $\g=-1$, which is extended to $-1<\g\leq 0$  \cite[Corollary 5.4]{GNROPUC}. 
Actually, since the LHS vanishes for $\nu = \GW(\g)$, we can also rewrite the sum rule \eqref{GWclassical} as
\begin{align}
\label{GWclassicalmodif}
\mathcal K(\GW(\g) \ | \ \nu ) & = \g \Re \sum_{k=0}^\infty (\alpha_k \bar\alpha_{k-1} - \alpha_k^\g \bar\alpha_{k-1}^\g)  - \sum_{k=0}^\infty \log\frac{1- |\alpha_k|^2}{1- |\alpha_k^\g|^2}\,,
\end{align}
where   $\alpha_k^\g$ is in \eqref{alphalimGW-}. This RHS is also the RHS of a sum rule for $\GW(\g)$ with $|\g| > 1$ (see \cite{GNRfutur}).
\end{remark}

\subsubsection{Hua-Pickrell}\label{susec:HP}

The Hua-Pickrell distribution appears in the study of random matrices corresponding to the potential
\begin{eqnarray}
\label{potHP2}
\V_\d(z) = -2\d \log |1-z|\,,
\end{eqnarray}
which is invariant by $z \mapsto \bar z$. Here, $\d$ is a complex parameter.
It has been introduced in \cite{hua1963harmonic} and appeared later in \cite{pickrell1987measures}. 
We also refer to \cite{Ner1}, 
 \cite{BO} and \cite{BNR}. 
 We will consider here only the case of real parameter $\d > 0$.

The equilibrium measure is the measure 
\begin{equation}
\label{limmeas}
\HP (\d) (dz) = (1+\d) \frac{\sqrt{\sin^2(\theta/2) -
\sin^2(\theta_\d/2)}}{2\pi \!\ \sin(\theta/2)} \mathbbm{1}_{(\theta_\d, 2\pi-\theta_\d)}(\theta) d\theta \,,
\end{equation} 
with $z= e^{\ii \theta}, \theta \in [0, 2\pi]$ and where 
$\theta_\d \in (0, \pi)$ is such that 
\begin{equation}\sin (\theta_\d/2) 
 = \frac{\d}
{ 1 + \d}\,. 
\end{equation}

The orthogonal polynomials with respect to $\HP (\d)$ are the Geronimus polynomials with constant V-coefficients
\begin{equation}
\label{defgammad}
\alpha_k \equiv \gamma_\d :=  - \frac{\d}{1+\d} \ , \ k \geq 0\,. 
\end{equation}
For $\gamma \in \mathbb{D}$, let 
\begin{align}\label{H1}
H_\d(\gamma) = -\log \frac{1 - |\gamma|^2}{1 - \gamma_\d^2}  - 2 \d \log \frac{|1- \gamma|}{1- \gamma_\d}\,.   
\end{align}

The arguments of the functions $H_\d$ in the sum rule are the deformed V-coefficients (see \cite[Section 2.2]{BNR}). For a nontrivial measure $\nu\in \mathcal M_1(\mathbb T)$ they form a sequence of parameters $\gamma_k\in \mathbb{D}$, $k\geq 0$ defined by
\begin{equation}
\label{phioverphi}
 \gamma_k = \bar\alpha_k \frac{\Phi_k^*(1)}{\Phi_k(1)},\quad (k\geq 0).
 \end{equation}
and can be computed via  the recursive definition
\begin{equation}
\label{mVc}
\gamma_0 = \bar\alpha_0 , \quad  \gamma_k = \bar\alpha_k \prod_{j=0}^{k-1} \frac{1 - \bar\gamma_j}{1 - \gamma_j}, \quad (k\geq 1).
\end{equation}
 Of course, when $\nu$ is symmetric, then $\Phi_k^*(1)=\Phi_k(1)$ and $\alpha_k$ is real, so that the deformed V-coefficients are the genuine V-coefficients.

Furthermore, define the function $\mathcal{F}_{\HP}$ for $\theta \notin (\theta_\d,  2\pi- \theta_d)$ :
\begin{align}
\label{defoutlierHP}
\mathcal F_{\HP} (e^{\ii \theta}) := \int_{I(\theta)}(1+\d) \frac{\sqrt{\sin^2\big(\theta_\d/2\big) -
\sin^2(\varphi/2)}}{\sin(\varphi/2)}\!\   d\varphi
\end{align}
where $I(\theta) = [\theta, \theta_\d]$ if $\theta \in (0, \theta_\d]$ and $I(\theta) = [2\pi-\theta_\d, \theta]$ if $\theta \in [2\pi-\theta_\d , 2\pi)$.

Then the following sum rule holds.

\begin{theorem}\cite[Theorem 5.1]{GNROPUC}
\label{sumruleHP}
Let $\d \geq 0$ and $\nu \in \mathcal M_1(\mathbb T)$ be nontrivial with  
$(\gamma_k)_{k \geq 0}\in {\mathbb{D}}^{\mathbb N}$ the sequence of its deformed V-coefficients. Then, if $\nu \in \Sr_1^{\mathbb T}(\theta_\d,2\pi-\theta_\d)$, 
\begin{align}
\label{SRHP}
 \mathcal K(\HP (\d) | \nu) + \sum_{\lambda \in E(\nu)} \mathcal F_{\HP}(\lambda) = \sum_{k=0}^\infty H_\d(\gamma_k)\,,
\end{align}
where both sides may be infinite simultaneously. If $\mu \notin \Sr_1^{\mathbb T}(\theta_\d,2\pi-\theta_\d)$, the RHS equals $+\infty$. 
\end{theorem}

\subsubsection{Poisson}

The Poisson kernel is the probability measure on $\mathbb T$ given by 
\begin{align}
\label{Pkern}
\Pois(\zeta) (dz) = \frac {1- |\zeta|^2}{2\pi|z-\zeta|^2} dz\,.
\end{align}
It is the equilibrium measure of random matrices with potential 
\begin{align}
\label{potP}
\mathcal V_\zeta (z) = \log| z\bar\zeta - 1|^2 ,
\end{align}
see \cite[Proposition 5.3.9]{HiaiP}, or \cite{hua1963harmonic} and \cite{baker1998random} for the study of the random matrix ensembles. Note that $\Pois(0)=\UNIF$.

The V-coefficients of $\Pois(\zeta)$ are
\begin{align} \label{PoissonVcoeff}
\alpha_0=  \zeta , \quad \alpha_k = 0 \ (k \geq 1)\,.
\end{align} 

We are aware of two sum rules relative to the Poisson measure $\Pois(\zeta)$.
The first one (Theorem 2.5.1 and formula (2.2.77) in \cite{simon05}) is
\begin{align*}
\mathcal K(\Pois(\zeta) | \nu) &= -\log\lambda_\infty (\zeta),
\end{align*}
where, with $\varphi_n$ the $n$-th orthonormal polynomial with respect to $\nu$, 
\begin{align*}
\lambda_\infty(\zeta) = (1-|\zeta|^2)\lim_{n\to\infty} |\varphi_n^*(\zeta)|^{-2}\,.
\end{align*}
The second one is quoted in Proposition \ref{rembess}. Its statement needs some notations given later in the paper.
We state a third new sum rule in Theorem \ref{Pois}.

\section{Mappings}
\label{sec:mappings}

Apart from the last one, all the mappings presented here  are from $\T$ to $\R$. They push forward a probability measure $\nu$ on the circle to a probability measure $\mu$ on the real line, which implies a possible connection of the J-coefficients of $\mu$ in terms of the V-coefficients of $\nu$. 
This may induce a connection between $J_\mu$ and $\mathcal C_\nu$ and also a correspondence, or "gateway" between  sum rules.

\subsection{Szeg{\H o}}
The Szeg{\H o} mapping from $\mathbb T$ to $[-2, 2]$ is defined by 
\begin{align}
\label{Smap}
 z &\mapsto \hbox{Sz}(z) = z+ z^{-1}\,,
\end{align}
or in angular coordinates,
\begin{align}
\label{Smapa}
\hbox{Sz}(e^{\ii\theta})= 2\cos\theta\,.
\end{align} 
The mapping $\hbox{Sz}$ is two-to-one from $\mathbb{T}$ to $[-2,2]$.  
For a symmetric $\nu \in \mathcal{M}_{1,s}(\mathbb{T})$, we let $\hbox{Sz}(\nu) = \nu \circ \Sz^{-1}$ be the pushforward of $\nu$ by the Szeg{\H o} mapping, which induces a bijection between $\mathcal{M}_{1,s}(\mathbb{T})$ and $\mathcal M_1([-2,2])$, the set of probability measures on $[-2,2]$. 
The OPUC $(\varphi_n)_{n\geq 0}$ with respect to $\nu$ and the OPRL $(p_n)_{n\geq 0}$ with respect to  $\hbox{Sz}(\nu)$ are related by
\begin{align}
\label{OPRLSz}
p_n(z) = z^{-n}\frac{\varphi_{2n}(z) + \varphi_{2n}^*(z)}
{\sqrt{2(1-\alpha_{2n-1})}},
\end{align}
where $\alpha_{2n-1}$ are the real V-coefficients of $\nu$. The Geronimus relations \cite[Theorem 13.1.7]{simon2}
  and equation (\ref{eq:coeffdecomposition})
give the remarkable identity
\begin{align} \label{Gero}
u_k(\Sz(\nu)) = \alpha_{k-1}(\nu)
\end{align}
for $k\geq 1$ between the canonical moments of $\Sz(\nu)$ and the V-coefficients of $\nu$.

\subsection{Delsarte-Genin (DG)}
For  ${\frak d}>0$ we consider the following relation between a point $z\in \mathbb{T}$ and $x\in [-2 {\frak d},2{\frak d}]$ given by    
\begin{align}
\label{DG0}
x =  \frak d\left(z^{1/2} + z^{-1/2}\right)\  \ \hbox{or} \ \ x = 2\frak d \cos(\theta/2)\,.
\end{align}
The following computations mainly come from \cite{delsartetri1}, \cite{delsartetri2}, \cite{delsarte1986split}. Therein, the parameter
$\frak d$ is fixed to $\tfrac{1}{2}$.
Notice that the  concern about branches of the square-root is addressed in  \cite[p. 518]{derevyagin2012cmv} .
With the right choice, this mapping is a bijection from $\mathbb{T}\setminus \{1\}$ to $(-2{\frak d},2{\frak d})$, which we denote by $\mathrm{DG}_{\frak d}$, the point $1\in \mathbb{T}$ corresponds to both $-2{\frak d}$ and $2{\frak d}$.  

Let $\nu\in \mathcal{M}_{s,1}(\mathbb{T})$ and 
fix ${\frak d}=1$. We let $\DG_1(\nu)$ be the pushforward of $\nu$ by $\DG_1$, with the convention that $\DG_1(\nu)(\{-2\})=\DG_1(\nu)(\{2\})=\tfrac{1}{2}\nu(\{1\})$. It is a symmetric measure on $[-2,2]$. 
The monic orthogonal polynomials $(\Phi_n)_{n\geq 0}$ with respect to $\nu$ and the monic orthogonal polynomials $(P_n)_{n\geq 0}$ with respect to  $\mu = \DG_1(\nu)$ are related by
\begin{align}
P_n(x) =\frac{z^{-n/2}(\Phi_n(z) + \Phi_n^*(z))}{\sqrt{2(1- \alpha_{n-1})}}\,,
\end{align}
where $\alpha_{n-1}$ are the real V-coefficients of $\nu$. The J-coefficients of $\mu$ are
\begin{align}
\label{a2alpha}
a_n^2 = (1+ \alpha_{n-1}) (1-\alpha_{n-2})  ,  \quad  b_n = 0 \qquad (n \geq 1).
\end{align}
The  canonical coefficients of $\mu$ of odd index are zero by symmetry, so comparing \eqref{eq:coeffdecomposition} and \eqref{a2alpha} we conclude
\begin{align}
\label{odd}u_{2n}(\DG_1(\nu)) =\alpha_{n-1}(\nu)  \  (n\geq 1) , \qquad u_{0}=\alpha_{-1}=-1 \,.
\end{align}
The inverse relation between $\Phi_n$ and $P_n$ is
\begin{align}
\Phi_n(z) = \frac{z^{n/2}\left(z^{1/2}P_{n+1}(x) - \sigma_n  P_n(x)\right)}{z-1},
\end{align}
with
\begin{align}
\label{n1}
\sigma_n =\frac{P_{n+1}(2)}{P_{n}(2)} = 1-\alpha_{n-1}.
\end{align}
An easy rescaling  is helpful when considering the general case ${\frak d}\neq 1$. Indeed, write the polynomials orthogonal to $ \DG_{\frak d}(\nu)$ as 
$P_n(x;\frak d) = \frak d^n P_n(x/\frak d)$. Then their V-coefficients satisfy 
\begin{align}
\label{alpha2a2}
a_n^2 = \frak d^2(1 + \alpha_{n-1})(1-\alpha_{n-2})\,.
\end{align}

Sometimes it is more convenient to use the mapping
\begin{align}
x =-i\frak d (z^{1/2} - z^{-1/2})\,,
\end{align}
or $x= 2\frak d \sin \theta$.
In this case we denote this mapping by $\DG_{\frak d}^-$ and the classical mapping by $\DG_{\frak d}^+$.

\subsection{Derevyagin-Vinet-Zhedanov (DVZ)}
\label{section DVZ}
This map was intruduced in \cite{derevyagin2012cmv} and generalized in \cite{cantero2021cmv}. It gives a remarkable relation between symmetric measures on $\mathbb{T}$ and measures on $\mathbb{R}$, induced by a algebraic relation between the CMV matrix and the Jacobi matrix, also called the Schur-Delsarte-Genin (SDG) map by \cite{cantero2021cmv}.

Let $\nu \in \mathcal{M}_{s,1}(\mathbb{T})$ be a symmetric measure. Its V-coefficients are real, and when the associated CMV matrix is written as  in \eqref{CMVdecomposition} in the form
$\mathcal C = \mathcal L\mathcal M$, 
we have in this case $\Theta_k^2=I_2$ (the identity in $\mathbb{R}^2$)  for all $k$. This implies
$\mathcal L^2 = \mathcal M^2 = I$ and the matrix 
$J_+ :=  \mathcal L + \mathcal M$ satisfies the following properties:
\begin{enumerate}
\item $J_+$ is real tridiagonal  symmetric.
\item $J_+^2 - 2 I = \mathcal C + \mathcal C^t$
\item The J-coefficients in $J_+$ are
\begin{align}
\label{abrho}
a_k =  \rho_{k-1} , \quad b_{k+1} = \alpha_{k} - \alpha_{k-1} \ \ (k \geq 1)
\end{align}
and $b_1=\alpha_0+1$. 
\item Its spectral measure is given by
\begin{align}
\label{defdvznu}
d\mu(x) = \frac{1}{2} (2 +x)\!\  d\DG_1(\nu) ,
\end{align}
supported on $[-2,2]$. Let us notice that this measure is not symmetric.
\end{enumerate}
The measure $\mu$ defined by (\ref{defdvznu}) will be denoted by $\DVZ^+ (\nu)$.

If we consider $J_- = \mathcal L - \mathcal M$ then the spectral measure satisfies
\begin{align}
\label{defdvznu-}
d\mu(x) = \frac{1}{2} (2 -x)\!\  d\DG_1(\nu)\,,
\end{align}
and it is denoted by $\DVZ^-(\nu)$.

\subsection{M\"obius}
The M\"obius transform $m_{z_0}$ for $z_0\in \mathbb D = \{ z\in \mathbb C : |z| < 1\}$ is defined by
\begin{align}\label{eq:mappingmoebius}
m_{z_0}(z) = \frac{z-z_0 }{1 - \bar z_0 z}\,.
\end{align}
It is an automorphism of  $\mathbb D$, sending $z_0$ to $0$, or of $\T$. Its inverse is $m_{-z_0}$.

\section{Pushing forward measures}

\subsection{Transformation of $\UNIF$}
From the definitions we see easily that
\begin{align}
\Sz (\UNIF) = \DG_1(\UNIF) = \Arcsine 
\end{align}
To compute the $\DVZ$ transform of $\UNIF$,  let us  introduce  the following notation. For $a<b$  the measure $\mathcal D(a,b)$ (resp. 
$\mathcal D (b,a)$)
is supported by $(a,b)$ with density
\begin{align}
\frac{2}{\pi(b-a)} \sqrt{\frac{x-a}{b-x}}\quad  \left(\hbox{resp.}\quad  
\frac{2}{\pi(b-a)}\sqrt{\frac{b-x}{x-a}}\right)\,.
\end{align}
These measures are affine pushforwards of   
the beta-distribution with parameter 
$\frac{1}{2}, -\frac{1}{2}$ (resp. $-\frac{1}{2},\frac{1}{2}$).  
The measure $\mathcal D(2, -2)$ is also a shift of the Marchenko-Pastur distribution, in the hard edge case.

The associated orthonormal polynomials are (up to an affine change) Chebyshev of the third type (resp. fourth type).
We then have
\begin{align*}
\DVZ^+ (\UNIF) = \mathcal D(-2, 2) \quad \text{ and }\quad \DVZ^- (\UNIF) = \mathcal D(2, -2) .
\end{align*}
The J-coefficients of  $\mathcal D(-2, 2)$ are by \eqref{abrho}
\begin{align}
a_k = 1 , \quad b_{k+1}=0\ (k \geq 1), 
\end{align}
and $ b_1 = 1$. 

\subsection{Transformation of $\GW$}
From the density 
 \eqref{GWMeq}, we deduce 
for $ \g\leq- 1$
 \begin{align*}
\Sz (\GW(\g))(dx) = \displaystyle   \frac{|\g|}{2\pi} \sqrt{\frac{4|\g|^{-1}-2-x}{x+2}} \mathbbm{1}_{(-2,4|\g|^{-1}-2)}(x)\ dx
\end{align*}
or in other words
\begin{align}
\Sz(\GW(\g)) =\mathcal D(4|\g|^{-1}-1 , -2)\,.
\end{align}
For $|\g| \leq 1$ let us notice that  the Gross-Witten density \eqref{GW-} is the mixture:
\begin{align}\label{linearcb}\GW(\g) = |\g| \GW(\epsilon(g)) + (1 - |\g|)\UNIF\,,
\end{align}
where $\epsilon(\g)$ is the sign of $\g$. 
Since the Szeg\H o mapping acts linearly on measures, 
 we obtain the complete picture
\begin{align}
\label{SzGW}
\Sz(\GW(\g)) = \mu_\g =: \begin{cases}
\mathcal D(-2+4|\g|^{-1}, -2)&\hbox{if}\ \g \leq -1 , \\
|\g|\mathcal D(2, -2) + (1- |\g|)\Arcsine&\hbox{if}\  -1 \leq \g \leq 0 , \\
\g\mathcal D(-2, 2) + (1- \g)\Arcsine&\hbox{if}\  0 \leq \g \leq 1 , \\
\mathcal D(2-4\g^{-1}, 2)&\hbox{if}\ \g \geq 1 .
\end{cases}
\end{align}

\subsubsection{With $\DG$  when  $|\g|\geq 1$} 

For $\g <-1$, 
the change of variable
\[x = 2\sqrt{|\g|}\cos (\theta/2)\]
gives
\begin{align}
\DG_{\sqrt{|\g|}}(\GW(\g))= \SC
\end{align}
When $\g > 1$, the change of variable 
\[x= 2\sqrt{\g }\sin (\theta/2)\]
gives also $\DG^-_{\sqrt{\g }}(\GW(\g))=\SC$.

\subsubsection{With $\DG$  when  $|\g|\leq 1$}
Starting from (\ref{linearcb}) and since $\DG$ is linear, we get for $-1 \leq \g< 0$
\begin{align}
\label{DG2W}
\DG_1(\GW(\g)) = 
 \rho_\g :=  
|\g| \SC + (1 - |\g|)\Arcsine . 
\end{align}
The corresponding canonical moments are 
\begin{align}
u_{2k} = \alpha_{k-1}^\g, \quad u_{2k+1} = 0 \quad (k \geq 0)\,,
\end{align}
where $\alpha_k^\g$ is in (\ref{alphalimGW-}). 

Similarly the application of $\DG_1^-$ leads,
for $0<\g\leq 1$, to the mixture
\begin{align} \label{DG2W2}
\DG_1^-(\GW(\g)) = \g  \SC + (1 -  \g )\Arcsine . 
\end{align}

\subsection{Transformation of $\HP$}
The
change of variable $x = 2 \cos \theta$ in \eqref{limmeas} gives
\begin{align*}
\frac{\Sz(\HP(\d))(dx)}{d x} = \frac{2(1+ \d)}{2\pi}\frac{\sqrt{x_\d-x}}{(2-x)\sqrt{2+x}}  = \frac{1+ \d}{\pi}\frac{\sqrt{(x_\d-x)(2+x)}}{(4-x^2)}\,.
\end{align*}
where
\begin{align}\label{defxd}
x_\d= \frac{2(1+2\d -\d^2)}{1+2\d+\d^2}\,.\end{align}
We conclude that
\begin{align}
\label{SzHP}
\Sz (\HP(\d)) = \KMK(2\d, 0)\,,
\end{align}
(recall that $\KMK(2\d, 0)$ is supported 
on $[-2,x_\d ]$).

 The change of variable $x =2 \cos (\theta/2)$   in \eqref{limmeas} gives for the density of $\DG_1(\HP(\d))$
\begin{align}
\notag 
(1+\d) \frac{\sqrt{\cos^2(\theta_\d/2) - \cos^2(\theta/2)}}{2\pi \sin^2(\theta/2)}= (1+\d)\frac{\sqrt{4\cos^2(\theta_\d/2) - x^2}}{\pi (4 -x^2)}\, ,
\end{align}
so that we conclude
\begin{align}
\DG_1(\HP(\d)) = \KMK(\d)\, .
\end{align} 
which is supported by $[-\hat x_\d,\hat x_\d ]$, with 
\begin{align} \label{defxdhat}
\hat x_\d = 2\frac{\sqrt{1+2\d}}{1+\d} .
\end{align}
Note that the V-coefficients of $\HP(\d)$ are constant equal to $\gamma_\d$ and then by \eqref{a2alpha} the J-coefficients of $\KMK(\d)$ are
\begin{align*}
 a_1^2 = 2 (1 + \gamma_\d) ,\ \  a_n^2 =  (1-\gamma_\d^2)\ \  (n \geq 2) ,\quad b_n = 0\ \ (n \geq 1)\, ,
\end{align*}
which agrees with \eqref{KMKcanonical}. 

Let us summarize the above results by two tables :   

\bigskip

\centerline{
\begin{tabular}{|c|c|}
\hline
&$\Sz$\\
\hline
$\UNIF$&$\Arcsine$\\
\hline 
$\GW(\g) , |\g|\leq 1$&$|\g|\mathcal D (-2\epsilon(\g), 2\epsilon(g)) + (1 - |\g|)\Arcsine$\\
\hline
$\GW(\g) , |\g| > 1$&$\mathcal D (2\epsilon(\g)- 4\g^{-1} ,2\epsilon(\g))$\\
\hline
$\HP(\d)$&$\KMK(2\d, 0)$\\
\hline
\end{tabular}}
\bigskip

\centerline{
\begin{tabular}{|c|c|c|}
\hline
&$\DG_1^+$&$\DG_1^-$\\
\hline
$\UNIF$&$\Arcsine$&$\Arcsine$\\
\hline 
$\GW(\g) , -1 \leq |\g|\leq 0$&$|\g|\SC + (1 - |\g|)\Arcsine$&\\
\hline
$\GW(\g) , 0 \leq |\g|\leq 1$&&$|\g|\SC + (1 - |\g|)\Arcsine$\\
\hline
$\GW(\g) , \g \leq -1 $&$\SC$ ($\ast$)&\\
\hline
$\GW(\g) , \g > 1$&&$\SC$ ($\ast$)\\
\hline
$\HP(\d)$&$\KMK(\d)$&\\
\hline
\end{tabular}}
\medskip

Here $(\ast)$ means that  $\DG^\pm_{\sqrt{|\g|}}$ is applied.

\section{Gateways}

In this section, we 
highlight connections between different sum rules 
arising when measures are transformed by the mappings of Section \ref{sec:mappings}. Unlike the large deviation technique developed in \cite{GaNaRo}, this method 
designing new
sum rules is purely analytical. 
Nevertheless, it
requires an existing sum rule to run.
In most cases, the aim is to obtain an OPRL sum rule from an OPUC one.

A measurable mapping $\varphi:X\to Y$ between metric spaces induces a mapping from $\mathcal{M}_1(X)$ to $\mathcal{M}_1(Y)$ by $\mu\mapsto \varphi(\mu)=\mu \circ \varphi^{-1}$. A sum rule for measures in $\mathcal{M}_1(Y)$ may lead to an identity for measures $\mu \in \mathcal{M}_1(X)$ (or vice versa) by evaluating both sides of the sum rule for $\varphi(\mu)$. 

Suppose the mapping $\varphi:X\to Y$ is a bijection. Then $\mu \mapsto \varphi(\mu)$ is a bijection from $\mathcal{M}_1(X)$ to $\mathcal{M}_1(Y)$. The entropy part of a sum rule can then
be obtained directly by the reversible entropy principle: 
\begin{align} \label{KLidentity}
\mathcal K(\mu_0 \mid \mu) = \mathcal K( \varphi(\mu_0) \mid  \varphi(\mu)) .
\end{align}
Among the mappings introduced in Section \ref{sec:mappings}, only the M\"obius mapping is one to one.
However, the Szeg\H{o} mapping is a bijection between symmetric measures on $\mathbb{T}$ and measures on $[-2,2]$ and \eqref{KLidentity} still holds for $\varphi=\Sz$ and $\mu,\mu_0\in \mathcal{M}_{1,s}(\mathbb{T})$. The mappings $\DG^\pm_{\frak d}$ are bijective when restricted to $\mathbb{T}\setminus \{1\}$. With the convention on mapping the mass at 1, \eqref{KLidentity} also holds for $\varphi =   \DG^\pm_{\frak d}$. Then \eqref{KLidentity}  
also holds
for $\varphi = \DVZ^\pm$ and therefore for all mappings considered in this paper.  

Transforming the RHS of a sum rule is less straightforward, but may be simplified if the coefficients  of $\varphi(\mu)$ are connected with those of $\mu$ in a convenient way.

\subsection{From $\UNIF$ to $\Arcsine$}
\label{susec:gatewayUNIFarcsine}

Using  \eqref{KLidentity} we get, with $\nu \in \mathcal{M}_{1,s}(\mathbb{T})$ such that $\Sz(\nu)=\mu$,
\begin{align*}
\mathcal K(\Arcsine \, | \, \mu)  = \mathcal K(\UNIF \, | \, \nu) .
\end{align*}
The Szeg{\H o} formula \eqref{SVsum},  jointly with \eqref{Gero}, gives
\begin{align*}
\mathcal K(\UNIF \, | \,  \nu) = -\sum_{k=0}^\infty \log (1 - |\alpha_k|^2(\nu))
=  -\sum_{k=1}^\infty \log (1 - u_k^2(\mu))\,,
\end{align*}
 and we recover the sum rule relative to the arcsine law \eqref{SRarcsine}.
Notice that $\UNIF$ is a particular case of the following distributions,
\[\UNIF = \HP(0) = \GW(0). \]
So that, the sum rule relative to $\Arcsine$ can be recovered from any sum rule relative to one of these distributions.

\subsection{From $\HP$ to $\KMK$}
\label{susec:gatewayHPKMK}

The Kesten-McKay laws can be obtained from the Hua-Pickrell distribution as 
\[\Sz(\HP(\d)) = \KMK(2\d, 0) , \quad \DG_1(\HP(\d)) = \KMK(\d).\]
These distributional identities allow to  
recover the sum rules with reference measure $\KMK(2\d, 0)$ or $\KMK(\d)$. 
If $\nu$ is a symmetric distribution in $\Sr_1^\mathbb T (\theta_\d, 2\pi-\theta_\d)$ 
then $\Sz(\HP(\d))$  (resp. $\DG_1(\nu)$) is supported on $[-2, 2]$  and belongs to $\Sr_1^\mathbb R (-2, x_\d)$ (resp. $\Sr_1^\mathbb R (-\hat x_\d, \hat x_\d))$. 

Let us consider the Szeg\H{o} mapping of the  sum rule \eqref{SRHP}. 
Let  $\mu$ be supported on $[-2, 2]$. If $\mu \in \Sr_1^\mathbb R (-2, x_\d)$, then a symmetric measure $\nu$ with $\Sz(\nu)=\mu$ belongs to $\Sr_1^\mathbb T (\theta_\d, 2\pi-\theta_\d)$. The LHS of the sum rule \eqref{SRHP} evaluated at $\nu$ is
\begin{align} \label{SRHPtransformedLHS}
  \mathcal K(\HP (\d) | \nu ) + \sum_{\lambda \in E(\nu)} \mathcal F_{\HP(\d)}(\lambda)
  = \mathcal K(\KMK(2\d, 0) | \mu ) + 2 \sum_{\lambda \in E(\mu)} \mathcal F(\lambda)
\end{align}
where the factor 2 comes from the two support points in $E(\nu)$ corresponding to one support point in $E(\mu)$, and where 
\begin{align}
2 \mathcal F (x) = 2 \mathcal F_{\HP(\d)}(e^{\ii \arccos (x/2)})=2 \int_{x_\d}^x  (1+\d)\frac{\sqrt{t- x_\d}}{(2-t)\sqrt{2+t}}\ dt\, = \mathcal{F}_{\KMK (2\d,0)}(x).
\end{align}
We therefore obtain (a particular case of) the LHS of \eqref{sumrule}. 
When looking at the RHS of \eqref{SRHP} evaluated at the symmetric measure $\nu$, we see that the deformed V-coefficients $\gamma_k$ are the regular V-coefficients $\alpha_k$. 
But by the 
 relation \eqref{Gero}, they are the canonical moments $u_k$ of $\mu=\Sz(\nu)$. Consequently, the RHS becomes
\begin{align}
\label{com1}
 - \sum_{k=1}^\infty \log\frac{1 - u_k^2}{1-\gamma_\d^2} - 
2\d \sum_{k=1}^\infty \log \frac{1-u_k}{1-\gamma_\d} 
 = -(1+2\d) \sum_{k=1}^\infty \log \frac{1-u_k}{1- \gamma_\d} -  \sum_{k=1}^\infty \log\frac{1 + u_k}{1+\gamma_\d}\, .  
\end{align}
The measure $\KMK(2\d,0)$ has all canonical moments equal to $\gamma_\d$, and so \eqref{com1} is equal to $I_J (\mu)$ as given in \eqref{rateJacobicoefficient}. We have recovered the sum rule corresponding to $\KMK(2\d, 0)$.

Let us consider the $\DG$ mapping and let $\mu$ supported on $[-2, 2]$ and symmetric. If $\nu$ is a symmetric measure on $\mathbb{T}$ such that $\DG_1(\nu) = \mu$ and $\nu \in  \Sr_1^\mathbb T (\theta_\d, 2\pi-\theta_\d)$, then $\mu \in  \Sr_1^\mathbb R (-\hat x_\d, \hat x_\d)$. The LHS of \eqref{SRHP} can then be transformed as 
\begin{align} \label{SRHPtransformedLHS2}
 & \mathcal K(\HP (\d) | \nu ) + \sum_{\lambda \in E(\nu)} \mathcal F_{\HP(\d)}(\lambda)\notag  \\
 & = \mathcal K(\KMK (\d) | \mu ) + \sum_{\lambda \in E(\mu)} \mathcal F(\lambda) .
\end{align}
where for $x>\hat x_\d$
\begin{align*}
\mathcal{F}(x) = \mathcal{F}_{\HP(\d)}(e^{2\ii \arccos (x/2)} ) = \mathcal{F}_{\KMK(\d)}(x) ,
\end{align*}
and we obtain the LHS of the $\KMK$ sum rule \eqref{sumrule}. Turning to the RHS, 
we first observe that the canonical moments of $\mu$ 
satisfy by \eqref{odd}
\[u_{2k+1}(\mu) =0 \ ,\ u_{2k}(\mu) =  \alpha_{k-1}(\nu)\,,\]
hence
\begin{align*}
H_\d(\gamma_{k-1}) &= H_\d(\alpha_{k-1}) =  -\log \frac{1-\alpha_{k-1}^2}{1-\gamma_\d^2} -2\d \log\frac{1-\alpha_{k-1}}{1-\gamma_\d}\\
&= -(1+2\d) \log\frac{1-u_{2k}}{1-\gamma_\d} - \log\frac{1+ u_{2k}}{1+\gamma_\d}  \\
&=\mathcal H_e^{\kappa, \kappa}(u_{2k})
\end{align*}
since $\gamma_\d  = u_e^{\d, \d}$. 
For odd index, we have $\mathcal H_o^{\d, \d}(u_{2k-1})=\mathcal H_o^{\d, \d}(0)=0$, since both $\mu$ and the reference measure are symmetric. 
We conclude that the RHS of the $\DG$ sum rule transforms exactly to the RHS of the $\KMK$ sum rule for symmetric measures.

\subsection{From $\GW$}

We now discuss
how 
one may obtain new sum rules starting from the sum rule \eqref{GWclassicalmodif} relative to $\GW(\g)$, 
for $-1\leq \g\leq 0$. Applying the Szeg\H{o} mapping leads to 
a sum rule relative to $\mu_\g$ in \eqref{SzGW}, a mixture of beta distributions. On the other hand, the mapping $\DG_1$ leads to a sum rule relative to a mixture of $\SC$ and $\Arcsine$. 

\subsubsection{With $\Sz$}
 
Let $\mu$ by a measure on $[-2,2]$ and $\nu \in \mathcal{M}_{1,s}(\mathbb{T})$ with $\Sz(\nu)=\mu$. Then the LHS of \eqref{GWclassical} applied to $\nu$ gives by \eqref{KLidentity}
\begin{align}
\mathcal K(\GW(\g)  |  \nu ) = \mathcal K(\mu_\g |  \mu ) .
\end{align}
In the RHS of \eqref{GWclassical} evaluated at $\nu$ only real V-coefficients appear.  By the 
relation \eqref{Gero} we obtain the sum rule for $\mu$ with support $[-2,2]$:
\begin{align} \label{newSRGW}
\mathcal K(\mu_\g |  \mu ) & = H(\g)   +\frac{\g }{2}-\frac{\g }{2} \sum_{k=1}^\infty (u_k - u_{k-1})^2  \\ \notag 
&\quad + \sum_{k=1}^\infty - \log (1 - u_k^2)  + \g u_k^2  ,
\end{align}
where $u_k$ are the canonical moments of $\mu$. 
\begin{remark}
When $\g = -1$, 
 (\ref{newSRGW}) becomes
\begin{align} \label{newSRGWu}
\mathcal K(\mu_{-1} |  \mu )  = \frac{1}{2} - \log 2 + 
\frac{1}{2}\sum_{k=0}^\infty (u_{k+1} - u_k)^2
+ \sum_{k=1}^\infty - \log (1 - u_k^2)  - u_k^2  ,
\end{align}
which  is  a version of  \cite[formula (2.8.6)]{simon05}.
 But 
\begin{align}
\label{identif}\mu_{-1}= \mathcal{D}(2,-2) = T(\MP_1)\,,\end{align} and
 where $T:\xi \mapsto \xi -2$, so that, by (\ref{KLidentity})
 \begin{align}
 \mathcal K(\mu_{-1} |  \mu ) = \mathcal K(\MP_1 | T^{-1}(\mu)) 
 \end{align}
 The RHS of the sum rule corresponding to $ \mathcal K(\MP_1 | T^{-1}(\mu))$ uses coefficients $(z_k)$ associated to $T^{-1}(\mu)$. To get an expression in terms of the $(u_k)$, we notice that 
if $a_k, b_k$ are the J-coefficients of $\mu$, the J-coefficients of $T^{-1}(\mu)$ are
\begin{align} \label{MPshifted1}
\tilde a_k = a_k  ,\quad  \tilde b_k = b_k + 2 .
\end{align}
Applying to $a_k$ and $b_k$ the decomposition into the canonical moment $u_k$ of $\mu$ according to \eqref{eq:coeffdecomposition}, and using  the parameters $(z_k)$ defined in \eqref{zcoeff} we obtain 
the relations
\begin{align} \label{MPshifted4}
z_1(T^{-1}(\mu))=2(1+u_1)(\mu), \quad z_{k}(T^{-1}(\mu)) =(1-u_{k-1}(\mu))(1+u_k(\mu))  , \quad (k \geq 1).
\end{align} 
Combining 
 the sum rule \eqref{SRL} relative to $\MP_1$ and \eqref{MPshifted4}, we obtain the identity
\begin{align} \label{newSRGW-1}
\mathcal K(\mathcal D(2,-2) | \mu) 
= \sum_{k=1}^\infty \big((1-u_{k-1})(1+u_k)-\log [(1-u_{k-1})(1+u_k)]-1 \big)
\end{align}
for $\mu$ with support $[-2,2]$.

Let us compare the RHS of (\ref{newSRGW-1}) and (\ref{newSRGWu}) by direct calculation.
Denote by $S_N$ the partial sum, up to $N$, of the RHS in (\ref{newSRGW-1}). 
Obviously, we may write
\begin{align}
\label{sumN1}
S_N &= \frac{1}{2} - \log 2 + u_N - \log (1+u_N) - \frac{u_N^2}{2}\\
\label{sumN2}
&+\frac{1}{2}\sum_{k=1}^N (u_k- u_{k-1})^2 + \sum_{k=1}^{N-1}(-\log (1-u_k^2) - u_k^2)
\end{align}
(recall that $u_0 = -1$). If the RHS of (\ref{newSRGW}) is finite, the gem (\ref{gem4}) 
warrants that
\begin{align}
\label{gemu}
\sum_{k=1}^\infty u_k^4 < \infty \ ,  \ \sum_{k=1}^\infty (u_{k+1}- u_k)^2 < \infty\ .
\end{align}
Hence, in particular $u_N \to 0$ and the two sums in (\ref{sumN2}) converge. We therefore recover (\ref{newSRGW}). Conversely, if one of the conditions in (\ref{gemu}) is not satisfied, we have $S_N \to \infty$ since the RHS of  (\ref{sumN1}) is bounded below by $-2\log2 -1$
\end{remark}

\subsubsection{With $\DG$} 

Recall that $\DG_1(\GW(\g))=\rho_\g$ as in \eqref{DG2W}. If $\mu$ is symmetric and supported on $[-2,2]$ and $\nu$ is symmetric on $\mathbb{T}$ such that $\DG_1(\nu)=\mu$, then from (\ref{KLidentity})
\begin{align}
\label{nsr3}
\mathcal K(\rho_\g | \mu) & = \mathcal K(\GW(\g) | \nu)
\end{align}
Now, in the sum rule (\ref{GWclassical}) the V-coefficients of $\nu$ are real and using 
\eqref{odd}, we may rewrite 
the last identity as
\begin{align} \label{newsr2b}
\mathcal K(\rho_\g | \mu) = 
 H(\g)   +\frac{\g }{2}-\frac{\g }{2} \sum_{k=1}^\infty (u_{2k} - u_{2k-2})^2 + \sum_{k=1}^\infty - \log (1 - u_{2k}^2)  + \g u_{2k}^2. 
\end{align}
Here, $u_k$ is the $k$-th canonical moment of $\mu$. 
The following theorem gives an alternative form of the RHS, obtained by combining the two sum rules relative to $\SC$ and $\Arcsine$.

\begin{theorem}
\label{newrhog}
For $\mu$ symmetric and supported on $[-2,2]$ and $-1\leq \g\leq 0$,
\begin{align}
\label{nsr2}
\mathcal K(\rho_\g  |  \mu)= C_\g +| \g| \sum_{k=1}^\infty( a_k^2 - 1 - \log a_k^2) + (1-|\g|) \sum_{k=1}^\infty \log a_k^2
\end{align}
where 
\begin{align}
\label{cste}
C_\g = -|\g| (1 - \log 2) + 1 - \sqrt{1-\g^2} + \log \frac{1+ \sqrt{1-\g^2}}{2}\,.
\end{align}
\end{theorem}

\proof 
From \eqref{DG2W}, $\rho_\g = |\g| \SC + (1 - |\g|)\Arcsine$, 
and applying Proposition \ref{prop:mixing} 
we thus obtain
\begin{align}
\notag
\mathcal K(\rho_\g\ | \ \mu)
&= |\g| \big(\mathcal K(\SC |  \mu) - \mathcal K(\SC |  \rho_\g)\big)\\
&\quad + (1-|\g|) \big(\mathcal K(\Arcsine |  \mu) - \mathcal K(\Arcsine |  \rho_\g)\big) 
\end{align}
From the Killip-Simon sum rule (Theorem \ref{sumruleg}) and from \eqref{Zlatos} we know $\mathcal K(\SC  |   \mu)$ and $\mathcal K(\Arcsine  |   \mu)$ respectively as  functions of the J-coefficients. This gives the coefficient dependent part of the RHS of \eqref{nsr2}. To compute the constant $C_\g$ we use  \eqref{KLidentity}, so that  
\begin{align*}
C_\g  &= -|\g|\mathcal K\left(\SC  |   \rho_\g \right)  -(1- |\g|)\mathcal K\left(\Arcsine  |   \rho_\g\right)\\
&= -|\g|\mathcal K(\GW(-1)  |   \GW(\g)) -(1- |\g|)\mathcal K(\UNIF  |   \GW(\g)) .
\end{align*}
The final value of $C_\g$  is then calculated with the help of   \cite[formula (7.5)]{GNROPUC}.
\qed

\begin{remark} 
When $\g\in (-1,1]$, we may use the alternative formulation \eqref{GWclassicalmodif} and compute the RHS using \eqref{a2alpha} and \eqref{odd}.  
\end{remark}

\subsubsection{With $\DVZ$} 
Let us restrict again to the case $\g = -1$.

Since $\DG_1(\GW(-1)) = \SC$, the measure $\hat\mu=\DVZ(\GW(-1))$ on $[-2,2]$ is by \eqref{defdvznu} 
\begin{align*}
d\hat\mu(x) = \frac{2+x}{2}d \SC(x) = \frac{1}{2\pi} (2+x)^{3/2} (2-x)^{1/2} dx\,.
\end{align*} 
We can then transform the sum rule \eqref{GWclassical} as follows. Assume that $\mu$ is a nontrivial measure supported by $[-2,2]$ and is such that there exists $\nu\in \mathcal M_{s,1}(\mathbb T)$ such that  $\mu = \DVZ(\nu)$. We have by (\ref{KLidentity})  
\begin{align} \label{KSSRvariant}
\mathcal{K}(\hat \mu | \mu) & = \mathcal{K}(\GW(-1)| \nu )
\end{align}
As before the V-coefficients of $\nu$ are real and related to the J-coefficients $a_k,b_k$ of $\mu$ by \eqref{abrho}. The sum rule relative to $\mathcal{K}(\GW(-1)| \nu )$ can be rewritten as
\begin{align}
\mathcal{K}(\hat \mu | \mu)
  = H(-1) -\frac{1 }{2}+\frac{1 }{2} \sum_{k=0}^\infty b_k^2 + \sum_{k=1}^\infty a_k^2-\log(a_k^2)-1 ,
\end{align}
 The RHS is therefore $\mathcal{I}_H(\mu)$ as in the Killip-Simon sum rule, (Theorem \ref{sumruleg} ), plus the negative constant $H(-1) -\frac{1 }{2}= 1/2 - \log 2$. 
 Notice that this does not mean that the RHS may be negative, actually this formula holds for $\mu$ in a restricted class.

Alternatively, the sum rule relative to $\hat \mu$ may be obtained directly from Theorem \ref{sumruleg}, since 
\begin{align*}
\mathcal{K}(\hat \mu | \mu) & = \mathcal K(\GW(-1)|\nu) = \mathcal K(\DG_1(\GW(-1)) | \DG (\nu))\\
&=\int  \log \frac{d \SC}{d \DG (\nu)} d\SC\\
&=\int  \log \frac{d \SC}{d\DVZ (\nu)} d\SC + \int \log \left(1 + \frac{x}{2}\right) d \SC(x)\\
&=\mathcal K(\SC | \DVZ (\nu)) + \int \log \left(1 + \frac{x}{2}\right) d \SC(x)
\end{align*}
and, as in \cite[Exercise 2.6.4]{agz}, 
\begin{align*}\int \log  \left(1 + \frac{x}{2}\right) \SC(dx)  = -\log 2 +\frac{1}{2}.
\end{align*}

\subsection{From $\UNIF$  to $\Pois$}
In this section we investigate sum rules relative to the measure $\Pois(\zeta)$ with $\zeta\in \mathbb{D}$ as given in \eqref{Pkern}. 
The (reverse) entropy with respect to the Poisson measure is called the Arov-Krein entropy (see \cite{roitberg2019arov}). The M\"obius transform $m_\zeta$ defined in \eqref{eq:mappingmoebius} maps $\Pois(\zeta)$ to the uniform measure $\UNIF$. 
Using \eqref{KLidentity} and the Szeg\H{o} sum rule \eqref{SVsum} shows that
\begin{align}
\label{sum0}
\mathcal{K}(\Pois(\zeta) | \nu) =\sum_{k=0}^\infty - \log (1 - |\alpha_{k}(m_{\zeta}(\nu))|^2)\,.
\end{align}
In the following, we analyze this sum rule and obtain alternative expressions for the RHS. 

Let us recall the connection between the V-coefficients and the Schur function of a measure  \cite[Chapter 1]{simon05}.
First, the Caratheodory function of a probability measure $\nu$ on $\mathbb T$ is defined as
\begin{align*}
F(z) = \int \frac{e^{\ii \theta} + z}{e^{\ii \theta} -z} d\nu(\theta) = 1 + 2z \int\frac{d\nu(\theta)}{e^{\ii \theta} -z}\,.
\end{align*}
It is analytic on $\mathbb{D}$. The Schur function is then defined from the Caratheodory function as
\begin{align}\label{defSch} f(z) = \frac{1}{z} \frac{F(z) -1}{F(z) +1} = \frac{1}{z} - \frac{2}{z(F(z) +1)}\end{align}
and conversely we have
\begin{align}
F(z) = \frac{1+zf(z)}{1-zf(z)} , \quad \frac{F(z) -1}{z} = \frac{2f(z)}{1-zf(z)}\,.
\end{align}
The V-coefficients can be obtained from $f$ by the classical Schur algorithm:
\begin{align}
S(g)(z)  = \frac{1}{z}\frac{g(z) - g(0)}{1 - \overline{g(0)}g(z)}\\
S^{[0]}(g) = g \ , \ S^{[k+1]}(g) =S\circ S^{[k]}\\
\alpha_k =  S^{[k]}f(0)\,.
\end{align}
To tackle the Poisson case, we use an extension of the Schur algorithm named the Nevanlinna-Pick algorithm that is defined as follows.

For $\rho \in \mathbb D \setminus\{0\}$ set, as in \cite{njaastad2009wall}
\begin{align}
S_\rho (g)  =\frac{ \omega(\rho)}{m_\rho} \frac{g - g(\rho)}{1- \overline{g(\rho)} g} \ ,  \  \omega(\rho)= -\frac{\rho}{|\rho|}\\
S_\rho^{[0]} (g) = g \ , \ S_\rho^{[k+1]} = S_\rho\circ S_\rho^{[k]} \ \ (k \geq 0)\,,
\end{align}
where we recall that $m_\rho$ is defined in (\ref{eq:mappingmoebius}).
To simplify we set $S_0 = S$.

We have then the following new sum rule, whose proof is postponed to Section \ref{sproof}.

\begin{theorem}
\label{Pois}
For $\nu$ a nontrivial measure on $\mathbb{T}$ with Schur function $f$, 
\begin{align}
\label{PSR}
\mathcal K(\Pois(\zeta) | \nu) = - \log (1 -|m_{\bar\zeta}\circ f(\zeta)|^2) + \sum_{k=1}^\infty - \log ( 1 -  |S_{\zeta}^{[k]}(f)(\zeta)|^2)\,.
\end{align} 
\end{theorem}

\begin{remark}
The prefactor $\omega$ introduced in Nevanlinna-Pick theory for technical reasons of infinite product convergence can be omitted here. Noticing that  $S_\zeta (\omega g) =  \omega S_\zeta(g)$ we can set
\[\hat S_\zeta = \omega^{-1} S_\zeta\]
and get recursively  
\[\hat S_\zeta^{[2k]} = S_\zeta^{[2k]}   \  ,  \hat S_\zeta^{[2k+1]} = \omega S_\zeta^{[2k+1]}\]
so that the sum in (\ref{PSR}) 
also holds for
$\hat S_\zeta^{[k]}$ instead of $S_\zeta^{[k]}$.
\end{remark}

There is another Poisson sum rule which is a direct consequence of a recent formula of  Bessonov \cite{bessonov}. We neither are able
to give any probabilistic interpretation nor to recover it from a pushforward of some other sum rule. The proof of this Poisson sum rule is also postponed to Section
\ref{sproof}.

\begin{proposition}
\label{rembess}
For $\nu$ a nontrivial measure on $\mathbb{T}$ with Schur function $f$, 
\begin{align} \label{SRBess}
\mathcal K(\Pois(\zeta) | \nu) = \log \frac{|1 -\zeta f(\zeta)|^2}{(1-|\zeta|^2)(1-|f(\zeta)|^2)}
+  \sum_{k= 1}^\infty \log \frac{1 -|\zeta S_0^{[k]}(f)(\zeta)|^2}{1 - |S_0^{[k]}(f)(\zeta)|^2}\,.
\end{align}
\end{proposition}
The above series has positive terms since 
\begin{align*}
|1 - \zeta f|^2 - (1 - |\zeta|^2)(1 - |f|^2) = |\bar \zeta - f|^2 \geq 0 ,
\end{align*}
and $1 - |\zeta f|^2  > 1 - |f|^2$, so that all the terms in the RHS of \eqref{SRBess} are positive. We conclude with the \emph{gems} corresponding to the above sum rules.

\begin{remark}
From (\ref{PSR}) and (\ref{SRBess}) we deduce that the Kullback-Leibler divergence $\mathcal K(\Pois(\zeta) | \nu)$ is finite if and only if
\[\sum_{k=1}^\infty |S_{\zeta}^{[k]}(f)(\zeta)|^2 < \infty , \]
or equivalently,
\[\sum_{k=1}^\infty |S_0^{[k]}(f)(\zeta)|^2 < \infty\,.\]
\end{remark}

\section{Proofs of Theorem \ref{Pois} and Prop. \ref{rembess}}
\label{sproof}
\subsection{Proof of Theorem \ref{Pois}}
The OPUC theory is in many ways an approximation theory.
The information carried by the V-coefficients is the same as the one
carried by the iterated Schur functions evaluated in $0$. This relies on 
the Schur function and its derivatives 
 at 
$0$. If we are interested in  the Schur  function and its derivatives at another point $\zeta\in \mathbb D$ we fall into the extension of the OPUC theory called Orthogonal Rational Functions (ORF) theory. Our main source for the following developments is 
\cite{njaastad2009wall} (see also \cite{bultheel1999orthogonal}). 

We start with  the sequence of rational functions
\[ 1, m_\zeta , (m_\zeta)^2, ... , (m_\zeta)^n , ...\]
and we apply the Gram-Schmidt orthonormalization in $L^2(\nu)$ to get
\[1, \varphi_1^o, \varphi_2^o, ...\]
we put the superscript  o to stress on the ORF aspect. 
 Let $(\Phi_n^o)$ be the corresponding monic ORF's.

Actually, we have
\[
\int \overline{\Phi_j^o (z)} \Phi_k^o(z) d\nu(z) = \kappa_k \delta_{jk}
\]
and if we set $z= m_{-\zeta}(\tau)$ we get
\[\int \overline{\Phi_j^o\circ m_{-\zeta}(\tau)}\Phi_k^o\circ m_{-\zeta}(\tau) d\nu^o(\tau) = \kappa_k \delta_{jk},,
\]
where $\nu^o = \nu \circ m_{-\zeta} = m_\zeta(\nu)$ is the pushforward of $\nu$ by $m_{\zeta}$. 
Of course the $\Phi_k^o\circ m_{-\zeta}$'s are the monic OPUC 
 with respect to $\nu^o$.

At the level of V-coefficients we have
\begin{align}
\Phi_k^o\circ m_{-\zeta} (0) = \Phi_k^o(\zeta)  , \quad \alpha_{k-1}(m_{\zeta}(\nu) )=-
\overline{\Phi_k^o(\zeta)} 
\,.
\end{align}

Now, let us study the relation between
the Schur functions. We write $F$ and $f$ for the Caratheodory and Schur function of $\nu$ and $F^o$ and $f^o$ for the functions of $\nu^o=m_\zeta(\nu)$. With $\tau = m_{-\zeta}(z)$, we have then
\begin{align}\notag
F^o(z) &= 1 + 2z \int\frac{d\nu^o(\theta)}{e^{\ii \theta} -z}= 1 + 2z \int\frac{1-\bar\zeta e^{\ii \theta}}{e^{\ii \theta}(1-\bar\zeta z)-\zeta -z} d\nu(\theta)\\
&= \frac{1-z\bar\zeta}{1+z\bar\zeta} + \frac{2z(1 - \tau\bar\zeta)}{(1 +z\bar\zeta)} \int \frac{d\nu(\theta)}{e^{\ii \theta} - \tau }\,.
\end{align}
Since 
\begin{align}
\int \frac{d\nu(\theta)}{e^{\ii \theta} - \tau 
} = \frac{F(\tau) -1}{2\tau}\,,
\end{align}
we obtain
the relation
\begin{align}
F^o(z) &= \frac{1-z\bar\zeta}{1+z\bar\zeta} + \frac{z(1 - \tau\bar\zeta)}{\tau(1 +z\bar\zeta)}(F(\tau) -1) .
\end{align}
By \eqref{defSch}, this implies
\begin{align}
f^o(z) &= 
\frac{-\bar\zeta + (1-\tau\bar\zeta) \frac{f(\tau)}{1-\tau f(\tau)}}
{1 + z(1-\tau\bar\zeta) \frac{f(\tau)}{1-\tau f(\tau)}}=  \frac{f(\tau)-\bar\zeta}{1+ f(\tau)(z-\tau z\bar\zeta - \tau)}
= \frac{f(\tau) - \bar\zeta}{1- \bar\zeta f(\tau)} 
\end{align}
or in other words
\begin{align}
f^o =m_{\bar\zeta}\circ f\circ m_{-\zeta} \,.
\end{align}
The first V-coefficient is then
\begin{align*}
\alpha^o_0 = f^o (0) = \left(m_{\bar\zeta}\circ f\right) (\zeta)=\frac{f(\zeta) -\bar\zeta}{1-\bar\zeta f(\zeta)} \,.
\end{align*}
To compute the higher order coefficients, let us begin with two auxiliary results. Observing that 
\begin{align}
\notag
S_\alpha(m_\beta \circ h)(z)  &= \frac{\omega(\alpha)}{m_\alpha(z)}\frac{m_\beta\circ h (z) -  m_\beta\circ h (\alpha)}
{1 - \overline{m_\beta\circ h (\alpha)} m_\beta\circ h (z)}= \frac{\omega(\alpha)}{m_\alpha(z)}\frac{\frac{h(z)-\beta}{1- \bar\beta h(z)} - \frac{h(\alpha)-\beta}{1- \bar\beta h(\alpha)}}{1 - \frac{\overline{h(\alpha)}-\bar\beta}{1- \beta \overline{h(\alpha)}}\frac{h(z)-\beta}{1- \bar\beta h(z)}}\\
&
 = \varepsilon (\alpha, \beta, h) \frac{\omega(\alpha)}{m_\alpha(z)}\frac{h(z) - h(\alpha)}{1- \overline{h(\alpha)} h(z)}
\end{align}
with
\[\varepsilon (\alpha, \beta, h) = \frac{1-\beta\overline{h(\alpha)}}{1 - \bar\beta h(\alpha)} \in \mathbb T . \]
So that, we have obtained
the first auxiliary result
\begin{align}
\label{Sm}
S_\alpha(m_\beta \circ h) = \varepsilon (\alpha, \beta, h) S_\alpha(h)\,.
\end{align}
The second one is the following
\begin{align}
\label{422}
S_0(h\circ m_\gamma)(z) = \frac{-1}{z}\frac{h\circ m_\gamma(z) - h(-\gamma)}{1 - \overline{h(-\gamma)}h\circ m_\gamma(z)}= 
-\overline{\omega(\gamma)}S_{-\gamma} (h) ( m_{\gamma}(z)) \,.
\end{align}
We have therefore
\begin{align}
\notag
S_0(f^o)&=
S(m_{\bar\zeta}\circ f\circ m_{-\zeta}) \sur{=}{(\ref{422})}-\overline{\omega(-\zeta)} \left[S_{\zeta}(m_{\bar{\zeta}}\circ f)\right]\circ m_{-\zeta}\\
\label{Sfo}
&\sur{=}{(\ref{Sm})} \varepsilon_1 S_{\zeta}(f) \circ m_{-\zeta} , 
\end{align}
with $
 \varepsilon_1 =\overline{\omega(\zeta)} \varepsilon (\zeta, \bar\zeta, f)\in \mathbb{T}$, 
hence
\begin{align}
\alpha^o_1 = S_0(f^o) (0) = \varepsilon_1 S_{\zeta}(f) (\zeta)\,.
\end{align}
This representation can be iterated. Assuming that
\begin{align}
\label{HypR}
S_0^{[k]} (f^o) = \varepsilon_k \left[S^{[k]}_{\zeta}(f)\right]\circ m_{-\zeta} ,\quad  \text{with } |\varepsilon_k| = 1\,,
\end{align}
we have
since $S_\alpha(\varepsilon h) = \varepsilon S_\alpha(h)$ when $\varepsilon \in \mathbb T$
\begin{align}
\notag
S_0^{[k+1]} (f^o) &= S_0\left[S_0^{[k]} (f^o)\right] =\varepsilon_kS\left[\left[S^{[k]}_{\zeta}(f)\right]\circ m_{-\zeta}\right]\\
&\sur{=}{(\ref{422})}-\overline{\omega(-\zeta)} \varepsilon_k  \left[ S_{\zeta}\left[S^{[k]}_{\zeta}(f)\right]\right]\circ m_{-\zeta} =\varepsilon_{k+1} \left[S^{[k+1]}_{\zeta}(f)\right]\circ m_{-\zeta} .
\end{align}
Inductively, \eqref{HypR} holds for every $k \geq 0$ and
\begin{align}
\alpha_k(m_\zeta (\nu)) = \alpha_k^o= S_0^{[k]}(f^0)(0) = \varepsilon_k S^{[k]}_{\zeta}(f)(\zeta)\,.
\end{align}
This finished the proof, since $|\varepsilon_k|=1$. \qed

\subsection{Proof of Proposition \ref{rembess}}
Let us recall the
Bessonov formula of \cite[Theorem 1]{bessonov}. 
Let $\nu$ be a probability measure on $\mathbb T$ with Lebesgue decomposition 
\[d\nu = g_\nu dz+d\nu_s\] (with respect to the uniform measure). The Bessonov formula is
\begin{align*}
\log \int \frac{d\Pois(\zeta)}{dz}\, d\nu -\int (\log g_\nu)\!\  d\Pois(\zeta)  = \sum_{k=0}^\infty \log \frac{1 -|\zeta f_k(\zeta)|^2}{1 - |f_k(\zeta)|^2}\,.
\end{align*}
Here, we set $f_k = S_0^{[k]} (f)$. It allows the following slight transformation. 
Since $\UNIF$ and $\Pois (\zeta)$ are mutually absolute continuous, it follows that  
$ g_\nu \frac{dz}{d\Pois(\zeta)}$ is the density of the absolutely continuous part of $\nu$ with respect to $\Pois(\zeta)$. Hence we get,
\begin{align*}
\notag -\int (\log g_\nu)\!\ d\Pois(\zeta) &= 
- \int (\log g_\nu)\!\  \frac{d z}{d\Pois (\zeta)} d\Pois(\zeta) - \int \log \frac{d\Pois(\zeta) }{dz}d\Pois(\zeta) \\
&= \mathcal K(\Pois(\zeta) | \nu) - \mathcal K(\Pois(\zeta) | \UNIF)\,. 
\end{align*} 
Besides, using the beginning of the proof of Lemma 1 in \cite{bessonov}, we  have
\begin{align}
\notag
 \int \frac{d\Pois(\zeta)}{dz}\!\ d\nu &= \frac{1 - |\zeta f(\zeta)|^2}{|1-\zeta f(\zeta)|^2} ,
\end{align}
Transforming the Kullback-Leibler distance according to \eqref{KLidentity} with the M\"obius mapping $m_\zeta$, using \eqref{SVsum} and \eqref{PoissonVcoeff}, we also have
\begin{align*}
\mathcal K(\Pois(\zeta) | \UNIF) =\mathcal K(\UNIF | \Pois(\zeta)) = -\log (1- |\zeta|^2) .
\end{align*} 
Consequently, we can write 
 \begin{align*}
\notag
\mathcal K(\Pois(\zeta)\ |\ \nu) &= -\log (1- |\zeta|^2) - \log \frac{1- |\zeta f(\zeta)|^2}{|1 - \zeta f(\zeta)|^2}+  \sum_{k=0}^\infty \log \frac{1 -|\zeta f_k(\zeta)|^2}{1 - |f_k(\zeta)|^2} \\
 &= \log \frac{|1 -\zeta f(\zeta)|^2}{(1-|\zeta|^2)(1-|f(\zeta)|^2)}
+  \sum_{k= 1}^\infty \log \frac{1 -|\zeta f_k(\zeta)|^2}{1 - |f_k(\zeta)|^2}\, ,
\end{align*}
which is the claimed sum rule. \qed

\section{Appendix}

\subsection{Analytical proof of a weak version of the HP sum rule}

Up to our knowledge, no analytical proof of the sum rule (\ref{sumruleHP}) is known. Nevertheless, we can express the coefficient side in terms of the limiting orthogonal polynomials and then use some limit theorems in the OP literature to try to recover the entropy of the spectral side, at least when there are no outliers.
\begin{proposition}
\label{weakprop}
 If the probability measure
$\mu= h \HP_\d$
is  such that there exists a polynomial $Q$ such that  $Qh$ and $Qh^{-1}$ are bounded on the arc $(\theta_\d, 2\pi-\theta_\d)$, then
\begin{align}
\mathcal K(\HP(\d)\ | \ \mu) = \sum_{k=0}^\infty H_\d(\gamma_k) < \infty\,.
\end{align}
\end{proposition}
It is a weaker form of  (\ref{sumruleHP}) since we impose stronger  conditions on $\mu$.

\proof
We will put a superscript $\d$ to all quantities relative to the reference mea\-sure. 

\emph{Step 1: Rewriting the coefficient side.}
The Szeg{\H o} recursion \eqref{recpolycirc} with \eqref{phioverphi} implies:
\begin{align}
\label{idb}
\Phi_n (1) = \prod_{k=0}^{n-1} (1 - \gamma_k)\,,
\end{align}so that 
\begin{align*}
\sum_{k=0}^{n-1} \log |1 - \gamma_j| = \log |\Phi_n (1)|
\,,
\end{align*}
and then
\begin{align}
\label{q1}
\sum_{k=0}^{n-1} \log \frac{|1 - \gamma_j|}{1- \gamma_\d} = 
 \log\frac {|\Phi_n (1)|
}{|\Phi^\d_n (1)|
}\,.\end{align}
If we go back to orthonormal polynomials, we have
\begin{align*}
\Phi_n (t) = \kappa_n^{-1} \varphi_n (t), \quad \kappa_n^{-2} = \prod_{k=0}^{n-1} (1 - |\alpha_k|^2)\end{align*}
and since $|\gamma_k| = |\alpha_k|$,
\begin{align}
\label{q2}
\sum_{k=0}^{n-1} - \log \frac{1- |\gamma_k|^2}{1- \gamma_\d^2} = 2 \log \frac{\kappa_n}{\kappa^\d_n}\,.
\end{align}
So that
\begin{align} 
\label{challenge}
\mathcal S_n := \sum_{k=0}^{n-1} H_\d(\gamma_k) = 2 \log \frac{\kappa_n}{\kappa^\d_n} -2\d\log\frac {|\Phi_n (1)|
}{|\Phi^\d_n (1)|
}\,.
\end{align}
Coming back to the normalized polynomials, we thus obtain,
\begin{align} \label{crucial}
\mathcal S_n = \sum_{k=0}^{n-1} H_\d(\gamma_k) = 2(1+ \d) \log \frac{\kappa_n}{\kappa^\d_n} -2\d\log\frac {|\varphi_n (1)|}{|\varphi^\d_n (1)|} .
\end{align}

\emph{Step 2: Computation of the limit.} 
We use an extension to measures supported on an arc of the classical Mat\'e-Nevai-Totik result on the full unit circle \cite[Theorem 9.4.1]{simon2}. The result for an arc is due to Bello Hernandez and Lopez Lagomasino \cite{Bello1998}. Theorem 2 therein shows, that if
\[\mu= h \HP_\d\]
is  such that there exists a polynomial $Q$ such that  $Qh$ and $Qh^{-1}$ are bounded on the arc $a=(\theta_\d, 2\pi-\theta_\d)$, then 
\begin{align*}
\lim_{n\to\infty} \frac{\varphi_n(\zeta)}{\varphi^\d_n(\zeta)}= D_a(h, \zeta)  ,\quad \mbox{ and }\quad
\lim_{n\to\infty} \frac{\kappa_n}{\kappa_n^\d} = D_a(h, \infty)\,,
\end{align*}
uniformly on compact subsets of $\bar{\mathbb C} \setminus a$. 
Here, the subscript \textit{a} stands for "the arc". 
To understand the limit, we need some more notations (well detailed in \cite[Section 2.2]{Bello2001} ).

Let 
\[\eta(\tau) = \tau+ \sqrt{\tau^2 -1}\]
(with root such that $|\eta(\tau)| > 1$) 
be the conformal mapping of $\bar{\mathbb C} \setminus [-1, 1]$ onto $\bar{\mathbb C} \setminus \{z: |z| \leq 1\}$ such that $\eta(\infty) = \infty$ and $\eta'(\infty) > 0$. Set
\[\c = \cot (\theta_\d/2)\,.\]
Let also
\begin{align*}
\nu(\zeta) = \eta\left(\frac{i}{\c}\frac{\zeta+1}{\zeta - 1}\right)\,,
\end{align*}
be the conformal mapping from $\bar{\mathbb C} \setminus a$ onto $\bar{\mathbb C} \setminus \{z: |z| \leq 1\}$.
In particular
\begin{align*} 
\nu(1) = \infty \ , \ \nu(\infty) = \eta(i/\c) = i \sqrt{1+2\d}
\,.
\end{align*}
Following  
\cite[formula (10)]{Bello2001}, or   \cite[Lemma 9]{Bello1998}, we have the indentity
\begin{align*}
D_a(h, \zeta) = \frac{D(h, \nu(\zeta))|D(h, \eta(i/\c))|}{D(h, \eta(i/\c))}\,,
\end{align*}
where $D$ is a variant of the famous Szeg{\H o} function:
\begin{align*}
D(h, z) &= \exp\left\{ \frac{1}{4\pi}\int_{0}^{2\pi}\log [ h(\tau) ] \frac{e^{i\theta} +z}{e^{i\theta} -z} d\theta \right\}
\end{align*}
with
\begin{align}
\label{masterchange}
\cot\frac{\tau}{2}&= \c \cos \theta\,.
\end{align}
This yields, respectively
\begin{align*}
|D_a (h, 1)|& = \exp \left\{ - \frac{1}{4\pi}\int_0^{2\pi}\log h(\tau) d\theta\right\} \,, \\
|D_a (h, \infty)|& =\exp \left\{\frac{1}{4\pi}\int_0^{2\pi}\log [h(\tau)]
\Re \frac{e^{i\theta} + i \sqrt{1+2\d}}{e^{i\theta} - i \sqrt{1+2\d}} d\theta\right\} \\
&= \exp \left\{- \frac{1}{4\pi}\int_0^{2\pi}\log [h(\tau)]\frac{\d}{1+\d- \sqrt{1+2\d}\sin\theta}d\theta \right\} .
\end{align*}
Going back to \eqref{crucial}, we see that the limit as 
$n\to \infty$ exists and is given by 
\begin{align}
\label{found}
\mathcal S_\infty  & = \lim_{n\to \infty} \mathcal{S}_n \notag
=  2(1+ \d) \log |D_a (h, \infty)| -2\d\log | D_a(h, 1)| \\
& = - \frac{1}{2\pi}\int_0^{2\pi}\log [h(\tau)]  \frac{\d \sqrt{1+2\d}\sin \theta}{1+\d- \sqrt{1+2\d} \sin \theta}
 \ d\theta \notag \\
 & = - \frac{1}{2\pi} \int_0^{2\pi}\log [h(\tau)]  \frac{\d \cos(\theta_\d/2)\sin \theta}{1- \cos(\theta_\d/2) \sin \theta}
 \ d\theta \,, 
\end{align}
where $\tau$ and $\theta$ are connected by (\ref{masterchange}).
Actually, splitting this integral in two parts and using $\sin (-\theta) = -\sin \theta$ leads to,
\begin{align*}
\notag
\mathcal S_\infty 
&=- \frac{1}{2\pi}\int_0^\pi \log [h(\tau)] \d \cos(\theta_\d/2) \left( \frac{\sin \theta}{1- \cos(\theta_\d/2) \sin \theta} - \frac{\sin \theta}{1+ \cos(\theta_\d/2) \sin \theta}\right) d\theta \\
&= - \frac{1}{2\pi}\int_0^\pi \log [h(\tau)]\frac{2 \d \cos^2(\theta_\d/2)\sin^2 \theta}{1 - \cos^2(\theta_\d /2) \sin^2 \theta}\, d\theta\,.
\end{align*}
Now, we have successively
\begin{align}
\notag\
d\theta = \left(2 \c \sin^2(\tau/2) \sin \theta\right)^{-1} d\tau\,,\\
\notag
\sin \theta = \frac{\sqrt{\sin^2 (\tau/2) - \sin^2 (\theta_\d/2)}}{\cos(\theta_\d/2)\sin(\tau/2)}\,,\\
\notag
1 - \cos^2(\theta_\d/2) \sin^2 \theta = \frac{\sin^2 (\theta_\d/2)}{\sin^2 (\tau/2)}\,,
\end{align}
so that, using the values of $\c$ and $\sin (\theta_\d/2)$:
\begin{align}
\notag
\mathcal S_\infty& = - \int_{\theta_\d}^{2\pi - \theta_\d} \log [h(\tau)] (1+\d)\frac{\sqrt{\sin^2 (\tau/2) - \sin^2 (\theta_\d/2)}}
{2 \pi \sin(\tau/2)}\  d\tau\\ \notag
&= - \int
 \log h(\tau) \HP(d\tau) = \mathcal K(\HP(\d)\mid \mu)\,.
\end{align}
This ends the proof.
\qedhere

\subsection{Kullback-Leibler distances for mixtures}

Suppose $\mu_1$ and $\mu_2$ are probability measures on some measurable space $S$. The following proposition is useful in the study of sum rules relative to a mixture of $\mu_1$ and $\mu_2$.

\begin{proposition} \label{prop:mixing}
Let 
$\tau_1,\tau_2> 0$ with $\tau_1+\tau_2=1$. Then, 
\begin{align} \label{afffinite}
\mathcal K(\mu_i \mid \tau_1 \mu_{1} + \tau_2 \mu_{2}) < \infty,\;\; (i=1,2).
\end{align} 
Moreover, for 
any probability measure $\mu$ on $S$,
\begin{align}
\notag
\mathcal K (\tau_1\mu_1 + \tau_2\mu_2 | \mu) & = \tau_1 \mathcal K(\mu_1 | \mu) +
 \tau_2 \mathcal K (\mu_2 | \mu)\\
 \label{affK} & \quad - \tau_1  \mathcal K(\mu_1 | \tau_1\mu_1 + \tau_2\mu_2)-  \tau_2  \mathcal K(\mu_2 | \tau_1\mu_1 + \tau_2\mu_2) ,
\end{align}
where both sides are simultaneously finite or infinite.
\end{proposition}
 
\proof
Since, for $i=1,2$, $\mu_{i} \ll \tau_1 \mu_{1} + \tau_2 \mu_{2}$ and
\begin{align*}
 \frac{d\mu_{i}}{d(\tau_1 \mu_{1} + \tau_2 \mu_{2})} \leq \frac{1}{\tau_i} \,,
\end{align*}
we obtain \eqref{afffinite}. 

For the proof of \eqref{affK} let us begin with a useful (but obvious)
remark. If $\nu_1$ and $\nu_2$ are two probability measures, 
such that $\nu_1\ll \nu_2$, then
\begin{align*}
\label{basic}
K(\nu_1 \ | \ \nu_2) < \infty \Longleftrightarrow \int \left|\log \frac{d\nu_1}{d\nu_2} \right| d\nu_1 < \infty .
\end{align*}
This follows from the inequality $u (\log u)_- \leq 1/e$ for $u >0$. 

\noindent Now, if $\mathcal K(\tau_1 \mu_{1} + \tau_2 \mu_{2} \ | \ \mu) <\infty$, then
\begin{align*}
\int \left|\log \frac{d(\tau_1 \mu_{1} + \tau_2 \mu_{2})}{d\mu}\right| d(\tau_1 \mu_{1} + \tau_2 \mu_{2}) < \infty
\end{align*}
hence for $i=1,2$
\begin{align*}
\int \left|\log \frac{d(\tau_1 \mu_{1} + \tau_2 \mu_{2})}{d\mu}\right| d\mu_{i} < \infty 
\end{align*}
and
\begin{align*}
\int \log \frac{d(\tau_1 \mu_{1} + \tau_2 \mu_{2})}{d\mu} d\mu_{i} < \infty 
\end{align*}
and eventually
\begin{align*}
\mathcal K(\tau_1 \mu_{1} + \tau_2 \mu_{2} \ | \ \mu)= \sum_{i=1}^2 \tau_i\int \log \frac{d(\tau_1 \mu_{1} + \tau_2 \mu_{2})}{d\mu} d\mu_{i}
\end{align*}
Adding $\sum_{i=1}^2 \tau_i  \mathcal K(\mu_i \ | \ \tau_1 \mu_{1} + \tau_2 \mu_{2})$ we get
\begin{align*}
\notag
\mathcal K(\tau_1 \mu_{1} + \tau_2 \mu_{2} \ | \ \mu) &+ \sum_{i=1}^2 \tau_i  \mathcal K(\mu_i \ | \ \tau_1 \mu_{1} + \tau_2 \mu_{2}) \\
= & \sum_{i=1}^2 \tau_i \int \log \left( \frac{d(\tau_1 \mu_{1} + \tau_2 \mu_{2})}{d\mu} \times \frac{d\mu_{i}}{d(\tau_1 \mu_{1} + \tau_2 \mu_{2})}\right)d\mu_i\\
= & \sum_{i=1}^2 \tau_i \mathcal K(\mu_{i} \ | \ \mu)\,.
\end{align*}
Conversely, if $\mathcal K(\mu_i \ | \ \mu)$ for $i=1,2$ are finite, then $\mathcal K(\tau_1 \mu_{1} + \tau_2 \mu_{2} \ | \ \mu) $ is finite by convexity.
\qed 
\\
\\
\subsection*{Acknowledgement}
Support from the ANR-3IA Artificial and Natural Intelligence
Toulouse Institute is gratefully acknowledged.

\bibliographystyle{plain}
\bibliography{preprint-GNR-gateway}
\end{document}